\documentstyle[amssymb,12pt]{amsart}

\input texdraw

\newcommand{\arccot}{\operatorname{arccot}}			
\newcommand{\HypFsimple}[2]{\sideset{_#1}{_#2}{F}}
\newcommand{\HypF}[5]{\sideset{_#1}{_#2}{F}\!\left[\matrix  
#3\\#4\endmatrix;{\displaystyle #5}\right]}
\def\po#1#2{(#1)_{#2}}
\newcommand{\floor}[1]{\left\lfloor#1\right\rfloor}

\newtheorem{thm}{Theorem}

\newtheorem{lem}[thm]{Lemma}

\newtheorem{conj}{Conjecture}

\numberwithin{equation}{section}

\catcode`\@=11
\font\tenln    = line10
\font\tenlnw   = linew10

\newskip\Einheit \Einheit=0.5cm
\newcount\xcoord \newcount\ycoord
\newdimen\xdim \newdimen\ydim \newdimen\PfadD@cke \newdimen\Pfadd@cke

\newcount\@tempcnta
\newcount\@tempcntb

\newdimen\@tempdima
\newdimen\@tempdimb

\newdimen\@wholewidth
\newdimen\@halfwidth

\newcount\@xarg
\newcount\@yarg
\newcount\@yyarg
\newbox\@linechar
\newbox\@tempboxa
\newdimen\@linelen
\newdimen\@clnwd
\newdimen\@clnht

\newif\if@negarg

\def\@whilenoop#1{}
\def\@whiledim#1\do #2{\ifdim #1\relax#2\@iwhiledim{#1\relax#2}\fi}
\def\@iwhiledim#1{\ifdim #1\let\@nextwhile=\@iwhiledim
        \else\let\@nextwhile=\@whilenoop\fi\@nextwhile{#1}}

\def\@whileswnoop#1\fi{}
\def\@whilesw#1\fi#2{#1#2\@iwhilesw{#1#2}\fi\fi}
\def\@iwhilesw#1\fi{#1\let\@nextwhile=\@iwhilesw
         \else\let\@nextwhile=\@whileswnoop\fi\@nextwhile{#1}\fi}

\def\thinlines{\let\@linefnt\tenln \let\@circlefnt\tencirc
  \@wholewidth\fontdimen8\tenln \@halfwidth .5\@wholewidth}
\def\thicklines{\let\@linefnt\tenlnw \let\@circlefnt\tencircw
  \@wholewidth\fontdimen8\tenlnw \@halfwidth .5\@wholewidth}
\thinlines

\def\PfadDicke#1{\PfadD@cke#1 \divide\PfadD@cke by2 \Pfadd@cke\PfadD@cke  
\multiply\PfadD@cke by2}
\long\def\LOOP#1\REPEAT{\def\BODY{#1}\ITERATE}
\def\ITERATE{\BODY \let\next\ITERATE \else\let\next\relax\fi \next}
\let\REPEAT=\fi
\def\Punkt{\hbox{\raise-2pt\hbox to0pt{\hss$\ssize\bullet$\hss}}}
\def\DuennPunkt(#1,#2){\unskip
  \raise#2 \Einheit\hbox to0pt{\hskip#1 \Einheit
          \raise-2.5pt\hbox to0pt{\hss$\bullet$\hss}\hss}}
\def\NormalPunkt(#1,#2){\unskip
  \raise#2 \Einheit\hbox to0pt{\hskip#1 \Einheit
          \raise-3pt\hbox to0pt{\hss\twelvepoint$\bullet$\hss}\hss}}
\def\DickPunkt(#1,#2){\unskip
  \raise#2 \Einheit\hbox to0pt{\hskip#1 \Einheit
          \raise-4pt\hbox to0pt{\hss\fourteenpoint$\bullet$\hss}\hss}}
\def\Kreis(#1,#2){\unskip
  \raise#2 \Einheit\hbox to0pt{\hskip#1 \Einheit
          \raise-4pt\hbox to0pt{\hss\fourteenpoint$\circ$\hss}\hss}}

\def\Line@(#1,#2)#3{\@xarg #1\relax \@yarg #2\relax
\@linelen=#3\Einheit
\ifnum\@xarg =0 \@vline
  \else \ifnum\@yarg =0 \@hline \else \@sline\fi
\fi}

\def\@sline{\ifnum\@xarg< 0 \@negargtrue \@xarg -\@xarg \@yyarg -\@yarg
  \else \@negargfalse \@yyarg \@yarg \fi
\ifnum \@yyarg >0 \@tempcnta\@yyarg \else \@tempcnta -\@yyarg \fi
\ifnum\@tempcnta>6 \@badlinearg\@tempcnta0 \fi
\ifnum\@xarg>6 \@badlinearg\@xarg 1 \fi
\setbox\@linechar\hbox{\@linefnt\@getlinechar(\@xarg,\@yyarg)}%
\ifnum \@yarg >0 \let\@upordown\raise \@clnht\z@
   \else\let\@upordown\lower \@clnht \ht\@linechar\fi
\@clnwd=\wd\@linechar
\if@negarg \hskip -\wd\@linechar \def\@tempa{\hskip -2\wd\@linechar}\else
     \let\@tempa\relax \fi
\@whiledim \@clnwd <\@linelen \do
  {\@upordown\@clnht\copy\@linechar
   \@tempa
   \advance\@clnht \ht\@linechar
   \advance\@clnwd \wd\@linechar}%
\advance\@clnht -\ht\@linechar
\advance\@clnwd -\wd\@linechar
\@tempdima\@linelen\advance\@tempdima -\@clnwd
\@tempdimb\@tempdima\advance\@tempdimb -\wd\@linechar
\if@negarg \hskip -\@tempdimb \else \hskip \@tempdimb \fi
\multiply\@tempdima \@m
\@tempcnta \@tempdima \@tempdima \wd\@linechar \divide\@tempcnta \@tempdima
\@tempdima \ht\@linechar \multiply\@tempdima \@tempcnta
\divide\@tempdima \@m
\advance\@clnht \@tempdima
\ifdim \@linelen <\wd\@linechar
   \hskip \wd\@linechar
  \else\@upordown\@clnht\copy\@linechar\fi}

\def\@hline{\ifnum \@xarg <0 \hskip -\@linelen \fi
\vrule height\Pfadd@cke width \@linelen depth\Pfadd@cke
\ifnum \@xarg <0 \hskip -\@linelen \fi}

\def\@getlinechar(#1,#2){\@tempcnta#1\relax\multiply\@tempcnta 8
\advance\@tempcnta -9 \ifnum #2>0 \advance\@tempcnta #2\relax\else
\advance\@tempcnta -#2\relax\advance\@tempcnta 64 \fi
\char\@tempcnta}

\def\Vektor(#1,#2)#3(#4,#5){\unskip\leavevmode
  \xcoord#4\relax \ycoord#5\relax
      \raise\ycoord \Einheit\hbox to0pt{\hskip\xcoord \Einheit
         \Vector@(#1,#2){#3}\hss}}

\def\Vector@(#1,#2)#3{\@xarg #1\relax \@yarg #2\relax
\@tempcnta \ifnum\@xarg<0 -\@xarg\else\@xarg\fi
\ifnum\@tempcnta<5\relax
\@linelen=#3\Einheit
\ifnum\@xarg =0 \@vvector
  \else \ifnum\@yarg =0 \@hvector \else \@svector\fi
\fi
\else\@badlinearg\fi}

\def\@hvector{\@hline\hbox to 0pt{\@linefnt
\ifnum \@xarg <0 \@getlarrow(1,0)\hss\else
    \hss\@getrarrow(1,0)\fi}}

\def\@vvector{\ifnum \@yarg <0 \@downvector \else \@upvector \fi}

\def\@svector{\@sline
\@tempcnta\@yarg \ifnum\@tempcnta <0 \@tempcnta=-\@tempcnta\fi
\ifnum\@tempcnta <5
  \hskip -\wd\@linechar
  \@upordown\@clnht \hbox{\@linefnt  \if@negarg
  \@getlarrow(\@xarg,\@yyarg) \else \@getrarrow(\@xarg,\@yyarg) \fi}%
\else\@badlinearg\fi}

\def\@upline{\hbox to \z@{\hskip -.5\Pfadd@cke \vrule width \Pfadd@cke
   height \@linelen depth \z@\hss}}

\def\@downline{\hbox to \z@{\hskip -.5\Pfadd@cke \vrule width \Pfadd@cke
   height \z@ depth \@linelen \hss}}

\def\@upvector{\@upline\setbox\@tempboxa\hbox{\@linefnt\char'66}\raise
     \@linelen \hbox to\z@{\lower \ht\@tempboxa\box\@tempboxa\hss}}

\def\@downvector{\@downline\lower \@linelen
      \hbox to \z@{\@linefnt\char'77\hss}}

\def\@getlarrow(#1,#2){\ifnum #2 =\z@ \@tempcnta='33\else
\@tempcnta=#1\relax\multiply\@tempcnta \sixt@@n \advance\@tempcnta
-9 \@tempcntb=#2\relax\multiply\@tempcntb \tw@
\ifnum \@tempcntb >0 \advance\@tempcnta \@tempcntb\relax
\else\advance\@tempcnta -\@tempcntb\advance\@tempcnta 64
\fi\fi\char\@tempcnta}

\def\@getrarrow(#1,#2){\@tempcntb=#2\relax
\ifnum\@tempcntb < 0 \@tempcntb=-\@tempcntb\relax\fi
\ifcase \@tempcntb\relax \@tempcnta='55 \or
\ifnum #1<3 \@tempcnta=#1\relax\multiply\@tempcnta
24 \advance\@tempcnta -6 \else \ifnum #1=3 \@tempcnta=49
\else\@tempcnta=58 \fi\fi\or
\ifnum #1<3 \@tempcnta=#1\relax\multiply\@tempcnta
24 \advance\@tempcnta -3 \else \@tempcnta=51\fi\or
\@tempcnta=#1\relax\multiply\@tempcnta
\sixt@@n \advance\@tempcnta -\tw@ \else
\@tempcnta=#1\relax\multiply\@tempcnta
\sixt@@n \advance\@tempcnta 7 \fi\ifnum #2<0 \advance\@tempcnta 64 \fi
\char\@tempcnta}

\def\Diagonale(#1,#2)#3{\unskip\leavevmode
  \xcoord#1\relax \ycoord#2\relax
      \raise\ycoord \Einheit\hbox to0pt{\hskip\xcoord \Einheit
         \Line@(1,1){#3}\hss}}
\def\AntiDiagonale(#1,#2)#3{\unskip\leavevmode
  \xcoord#1\relax \ycoord#2\relax 
      \raise\ycoord \Einheit\hbox to0pt{\hskip\xcoord \Einheit
         \Line@(1,-1){#3}\hss}}
\def\Pfad(#1,#2),#3\endPfad{\unskip\leavevmode
  \xcoord#1 \ycoord#2 \thicklines\ZeichnePfad#3\endPfad\thinlines}
\def\ZeichnePfad#1{\ifx#1\endPfad\let\next\relax
  \else\let\next\ZeichnePfad
    \ifnum#1=1
      \raise\ycoord \Einheit\hbox to0pt{\hskip\xcoord \Einheit
         \vrule height\Pfadd@cke width1 \Einheit depth\Pfadd@cke\hss}%
      \advance\xcoord by 1
    \else\ifnum#1=2
      \raise\ycoord \Einheit\hbox to0pt{\hskip\xcoord \Einheit
        \hbox{\hskip-\PfadD@cke\vrule height1 \Einheit width\PfadD@cke  
depth0pt}\hss}%
      \advance\ycoord by 1
    \else\ifnum#1=3
      \raise\ycoord \Einheit\hbox to0pt{\hskip\xcoord \Einheit
         \Line@(1,1){1}\hss}
      \advance\xcoord by 1
      \advance\ycoord by 1
    \else\ifnum#1=4
      \raise\ycoord \Einheit\hbox to0pt{\hskip\xcoord \Einheit
         \Line@(1,-1){1}\hss}
      \advance\xcoord by 1
      \advance\ycoord by -1
    \fi\fi\fi\fi
  \fi\next}
\def\hSSchritt{\leavevmode\raise-.4pt\hbox to0pt{\hss.\hss}\hskip.2\Einheit
  \raise-.4pt\hbox to0pt{\hss.\hss}\hskip.2\Einheit
  \raise-.4pt\hbox to0pt{\hss.\hss}\hskip.2\Einheit
  \raise-.4pt\hbox to0pt{\hss.\hss}\hskip.2\Einheit
  \raise-.4pt\hbox to0pt{\hss.\hss}\hskip.2\Einheit}
\def\vSSchritt{\vbox{\baselineskip.2\Einheit\lineskiplimit0pt
\hbox{.}\hbox{.}\hbox{.}\hbox{.}\hbox{.}}}
\def\DSSchritt{\leavevmode\raise-.4pt\hbox to0pt{%
  \hbox to0pt{\hss.\hss}\hskip.2\Einheit
  \raise.2\Einheit\hbox to0pt{\hss.\hss}\hskip.2\Einheit
  \raise.4\Einheit\hbox to0pt{\hss.\hss}\hskip.2\Einheit
  \raise.6\Einheit\hbox to0pt{\hss.\hss}\hskip.2\Einheit
  \raise.8\Einheit\hbox to0pt{\hss.\hss}\hss}}
\def\dSSchritt{\leavevmode\raise-.4pt\hbox to0pt{%
  \hbox to0pt{\hss.\hss}\hskip.2\Einheit
  \raise-.2\Einheit\hbox to0pt{\hss.\hss}\hskip.2\Einheit
  \raise-.4\Einheit\hbox to0pt{\hss.\hss}\hskip.2\Einheit
  \raise-.6\Einheit\hbox to0pt{\hss.\hss}\hskip.2\Einheit
  \raise-.8\Einheit\hbox to0pt{\hss.\hss}\hss}}
\def\SPfad(#1,#2),#3\endSPfad{\unskip\leavevmode
  \xcoord#1 \ycoord#2 \ZeichneSPfad#3\endSPfad}
\def\ZeichneSPfad#1{\ifx#1\endSPfad\let\next\relax
  \else\let\next\ZeichneSPfad
    \ifnum#1=1
      \raise\ycoord \Einheit\hbox to0pt{\hskip\xcoord \Einheit
         \hSSchritt\hss}%
      \advance\xcoord by 1
    \else\ifnum#1=2
      \raise\ycoord \Einheit\hbox to0pt{\hskip\xcoord \Einheit
        \hbox{\hskip-2pt \vSSchritt}\hss}%
      \advance\ycoord by 1
    \else\ifnum#1=3
      \raise\ycoord \Einheit\hbox to0pt{\hskip\xcoord \Einheit
         \DSSchritt\hss}
      \advance\xcoord by 1
      \advance\ycoord by 1
    \else\ifnum#1=4
      \raise\ycoord \Einheit\hbox to0pt{\hskip\xcoord \Einheit
         \dSSchritt\hss}
      \advance\xcoord by 1
      \advance\ycoord by -1
    \fi\fi\fi\fi
  \fi\next}
\def\Koordinatenachsen(#1,#2){\unskip
 \hbox to0pt{\hskip-.5pt\vrule height#2 \Einheit width.5pt depth1 \Einheit}%
 \hbox to0pt{\hskip-1 \Einheit \xcoord#1 \advance\xcoord by1
    \vrule height0.25pt width\xcoord \Einheit depth0.25pt\hss}}
\def\Koordinatenachsen(#1,#2)(#3,#4){\unskip
 \hbox to0pt{\hskip-.5pt \ycoord-#4 \advance\ycoord by1
    \vrule height#2 \Einheit width.5pt depth\ycoord \Einheit}%
 \hbox to0pt{\hskip-1 \Einheit \hskip#3\Einheit
    \xcoord#1 \advance\xcoord by1 \advance\xcoord by-#3
    \vrule height0.25pt width\xcoord \Einheit depth0.25pt\hss}}
\def\Gitter(#1,#2){\unskip \xcoord0 \ycoord0 \leavevmode
  \LOOP\ifnum\ycoord<#2
    \loop\ifnum\xcoord<#1
      \raise\ycoord \Einheit\hbox to0pt{\hskip\xcoord \Einheit\Punkt\hss}%
      \advance\xcoord by1
    \repeat
    \xcoord0
    \advance\ycoord by1
  \REPEAT}
\def\Gitter(#1,#2)(#3,#4){\unskip \xcoord#3 \ycoord#4 \leavevmode
  \LOOP\ifnum\ycoord<#2
    \loop\ifnum\xcoord<#1
      \raise\ycoord \Einheit\hbox to0pt{\hskip\xcoord \Einheit\Punkt\hss}%
      \advance\xcoord by1
    \repeat
    \xcoord#3
    \advance\ycoord by1
  \REPEAT}
\def\Label#1#2(#3,#4){\unskip \xdim#3 \Einheit \ydim#4 \Einheit
  \def\lo{\advance\xdim by-.5 \Einheit \advance\ydim by.5 \Einheit}%
  \def\llo{\advance\xdim by-.25cm \advance\ydim by.5 \Einheit}%
  \def\loo{\advance\xdim by-.5 \Einheit \advance\ydim by.25cm}%
  \def\o{\advance\ydim by.25cm}%
  \def\ro{\advance\xdim by.5 \Einheit \advance\ydim by.5 \Einheit}%
  \def\rro{\advance\xdim by.25cm \advance\ydim by.5 \Einheit}%
  \def\roo{\advance\xdim by.5 \Einheit \advance\ydim by.25cm}%
  \def\l{\advance\xdim by-.30cm}%
  \def\r{\advance\xdim by.30cm}%
  \def\lu{\advance\xdim by-.5 \Einheit \advance\ydim by-.6 \Einheit}%
  \def\llu{\advance\xdim by-.25cm \advance\ydim by-.6 \Einheit}%
  \def\luu{\advance\xdim by-.5 \Einheit \advance\ydim by-.30cm}%
  \def\u{\advance\ydim by-.30cm}%
  \def\ru{\advance\xdim by.5 \Einheit \advance\ydim by-.6 \Einheit}%
  \def\rru{\advance\xdim by.25cm \advance\ydim by-.6 \Einheit}%
  \def\ruu{\advance\xdim by.5 \Einheit \advance\ydim by-.30cm}%
  #1\raise\ydim\hbox to0pt{\hskip\xdim
     \vbox to0pt{\vss\hbox to0pt{\hss$#2$\hss}\vss}\hss}%
}
\catcode`\@=13

\def\ldreieck{\bsegment
  \rlvec(0.866025403784439 .5) \rlvec(0 -1)
  \rlvec(-0.866025403784439 .5)
  \savepos(0.866025403784439 -.5)(*ex *ey)
        \esegment
  \move(*ex *ey)
        }
\def\rdreieck{\bsegment
  \rlvec(0.866025403784439 -.5) \rlvec(-0.866025403784439 -.5)  \rlvec(0 1)
  \savepos(0 -1)(*ex *ey)
        \esegment
  \move(*ex *ey)
        }
\def\rhombus{\bsegment
  \rlvec(0.866025403784439 .5) \rlvec(0.866025403784439 -.5)
  \rlvec(-0.866025403784439 -.5)  \rlvec(0 1)
  \rmove(0 -1)  \rlvec(-0.866025403784439 .5)
  \savepos(0.866025403784439 -.5)(*ex *ey)
        \esegment
  \move(*ex *ey)
        }
\def\RhombusA{\bsegment
  \rlvec(0.866025403784439 .5) \rlvec(0.866025403784439 -.5)
  \rlvec(-0.866025403784439 -.5) \rlvec(-0.866025403784439 .5)
  \savepos(0.866025403784439 -.5)(*ex *ey)
        \esegment
  \move(*ex *ey)
        }
\def\RhombusB{\bsegment
  \rlvec(0.866025403784439 .5) \rlvec(0 -1)
  \rlvec(-0.866025403784439 -.5) \rlvec(0 1)
  \savepos(0 -1)(*ex *ey)
        \esegment
  \move(*ex *ey)
        }
\def\RhombusC{\bsegment
  \rlvec(0.866025403784439 -.5) \rlvec(0 -1)
  \rlvec(-0.866025403784439 .5) \rlvec(0 1)
  \savepos(0.866025403784439 -.5)(*ex *ey)
        \esegment
  \move(*ex *ey)
        }
\def\hdSchritt{\bsegment
  \lpatt(.05 .13)
  \rlvec(0.866025403784439 -.5)
  \savepos(0.866025403784439 -.5)(*ex *ey)
        \esegment
  \move(*ex *ey)
        }
\def\vdSchritt{\bsegment
  \lpatt(.05 .13)
  \rlvec(0 -1)
  \savepos(0 -1)(*ex *ey)
        \esegment
  \move(*ex *ey)
        }

\begin{document}



\title[The number of rhombus tilings]{The number of rhombus
tilings of a symmetric hexagon which contain a fixed rhombus on the
symmetry axis, I}

\author{M.~Fulmek and C.~Krattenthaler}
\address{
Institut f\"ur Mathematik der Universit\"at Wien\\
Strudlhofgasse 4, A-1090 Wien, Austria\\
e-mail: Mfulmek@@Mat.Univie.Ac.At, Kratt@@Pap.Univie.Ac.At
}

\date{\today}

\begin{abstract}
We compute the number
of rhombus tilings of a hexagon with sides $N,M,N,\break N,M,N$, which
contain a fixed rhombus on the symmetry axis that cuts through the
sides of length $M$. 
\end{abstract}

\maketitle

\section{Introduction}

Let $a$, $b$ and $c$ be positive integers, and consider a 
hexagon of sides $a,b,c,a,b,c$ whose angles are $120^\circ$
(see Figure~\ref{fig:0}.a). The subject of our interest
is tilings of this hexagon
by rhombi of unit edge-length and angles of $60^\circ$ and
$120^\circ$ (see Figure~\ref{fig:0}.b). (From now on, by a
rhombus we always mean such a rhombus with side lengths 1 and angles of
$60^\circ$ and $120^\circ$.) 
By a well-known bijection \cite{DT}, the total number of rhombus tilings
of this hexagon is equal
to the number of all plane
partitions contained in an $a\times b\times c$ box. The
latter enumeration was solved long ago by MacMahon 
\cite[Sec.~429, $q\rightarrow 1$; proof in Sec.~494]{MacMahon}. Therefore:

\smallskip
{\it The number of all rhombus tilings of a hexagon
with sides $a,b,c,a,b,c$ equals}
\begin{equation}
\label{eq:MacMahon}
\prod_{i=1}^a\prod_{j=1}^b\prod_{k=1}^c\frac{i+j+k-1}{i+j+k-2}.
\end{equation}
(The form of the expression is due to Macdonald.)
\smallskip

\begin{figure}
\centertexdraw{
  \drawdim truecm  \linewd.02
  \rhombus \rhombus \rhombus \rhombus \ldreieck
  \move (-0.866025403784439 -.5)
  \rhombus \rhombus \rhombus \rhombus \rhombus \ldreieck
  \move (-1.732050807568877 -1)
  \rhombus \rhombus \rhombus \rhombus \rhombus \rhombus \ldreieck
  \move (-1.732050807568877 -1)
  \rdreieck
  \rhombus \rhombus \rhombus \rhombus \rhombus \rhombus \ldreieck
  \move (-1.732050807568877 -2)
  \rdreieck
  \rhombus \rhombus \rhombus \rhombus \rhombus \rhombus \ldreieck
  \move (-1.732050807568877 -3)
  \rdreieck
  \rhombus \rhombus \rhombus \rhombus \rhombus \rhombus 
  \move (-1.732050807568877 -4)
  \rdreieck
  \rhombus \rhombus \rhombus \rhombus \rhombus 
  \move (-1.732050807568877 -5)
  \rdreieck
  \rhombus \rhombus \rhombus \rhombus 
\move(8 0)
\bsegment
  \drawdim truecm  \linewd.02
  \rhombus \rhombus \rhombus \rhombus \ldreieck
  \move (-0.866025403784439 -.5)
  \rhombus \rhombus \rhombus \rhombus \rhombus \ldreieck
  \move (-1.732050807568877 -1)
  \rhombus \rhombus \rhombus \rhombus \rhombus \rhombus \ldreieck
  \move (-1.732050807568877 -1)
  \rdreieck
  \rhombus \rhombus \rhombus \rhombus \rhombus \rhombus \ldreieck
  \move (-1.732050807568877 -2)
  \rdreieck
  \rhombus \rhombus \rhombus \rhombus \rhombus \rhombus \ldreieck
  \move (-1.732050807568877 -3)
  \rdreieck
  \rhombus \rhombus \rhombus \rhombus \rhombus \rhombus 
  \move (-1.732050807568877 -4)
  \rdreieck
  \rhombus \rhombus \rhombus \rhombus \rhombus 
  \move (-1.732050807568877 -5)
  \rdreieck
  \rhombus \rhombus \rhombus \rhombus 
  \linewd.08
  \move(0 0)
  \RhombusA \RhombusB \RhombusB 
  \RhombusA \RhombusA \RhombusB \RhombusA \RhombusB \RhombusB
  \move (-0.866025403784439 -.5)
  \RhombusA \RhombusB \RhombusB \RhombusB \RhombusB
  \RhombusA \RhombusA \RhombusB \RhombusA 
  \move (-1.732050807568877 -1)
  \RhombusB \RhombusB \RhombusA \RhombusB \RhombusB \RhombusA
  \RhombusB \RhombusA \RhombusA 
  \move (1.732050807568877 0)
  \RhombusC \RhombusC \RhombusC 
  \move (1.732050807568877 -1)
  \RhombusC \RhombusC \RhombusC 
  \move (3.464101615137755 -3)
  \RhombusC 
  \move (-0.866025403784439 -.5)
  \RhombusC
  \move (-0.866025403784439 -1.5)
  \RhombusC
  \move (0.866025403784439 -2.5)
  \RhombusC \RhombusC 
  \move (0.866025403784439 -3.5)
  \RhombusC \RhombusC \RhombusC 
  \move (2.598076211353316 -5.5)
  \RhombusC 
  \move (0.866025403784439 -5.5)
  \RhombusC 
  \move (-1.732050807568877 -3)
  \RhombusC 
  \move (-1.732050807568877 -4)
  \RhombusC 
  \move (-1.732050807568877 -5)
  \RhombusC \RhombusC 
\esegment
\htext (-1.5 -9){\small a. A hexagon with sides $a,b,c,a,b,c$,}
\htext (-1.5 -9.5){\small \hphantom{a. }where $a=3$, $b=4$, $c=5$}
\htext (6.8 -9){\small b. A rhombus tiling of a hexagon}
\htext (6.8 -9.5){\small \hphantom{b. }with sides $a,b,c,a,b,c$}
\rtext td:0 (4.3 -4){$\sideset {}  c 
    {\left.\vbox{\vskip2.2cm}\right\}}$}
\rtext td:60 (2.6 -.5){$\sideset {} {} 
    {\left.\vbox{\vskip1.7cm}\right\}}$}
\rtext td:120 (-.44 -.25){$\sideset {}  {}  
    {\left.\vbox{\vskip1.3cm}\right\}}$}
\rtext td:0 (-2.3 -3.6){$\sideset {c}  {} 
    {\left\{\vbox{\vskip2.2cm}\right.}$}
\rtext td:240 (0 -7){$\sideset {}  {}  
    {\left.\vbox{\vskip1.7cm}\right\}}$}
\rtext td:300 (3.03 -7.25){$\sideset {}  {}  
    {\left.\vbox{\vskip1.4cm}\right\}}$}
\htext (-.9 0.1){$a$}
\htext (2.8 -.1){$b$}
\htext (3.2 -7.8){$a$}
\htext (-0.3 -7.65){$b$}
}
\caption{}
\label{fig:0}
\end{figure}

A natural question to be asked is what the distribution of the
rhombi in a random tiling is. 
On an {\it asymptotic} level, this question was answered by Cohn, 
Larsen and Propp \cite{CoLPAA}. On the {\it exact\/} (enumerative) level, 
Ciucu and
Krattenthaler \cite{CiucKratAB} and, independently, 
Helfgott and Gessel \cite{HeGeAA} 
computed the number of all rhombus
tilings of a hexagon with sides $N,M,N,N,M,N$ which contain the
central rhombus. (They were 
motivated by a problem posed by Propp \cite[Problem~1]{PropAA}.)

In this paper we solve an even more general problem,
namely the enumeration of all rhombus tilings of a hexagon with side lengths
$N,M,N,N,M,N$ which contain an {\it arbitrary} fixed rhombus
on the symmetry axis which cuts through the sides of length $M$ (see
Figure~\ref{fig:1} for illustration, the fixed rhombus is shaded). Our
results are the following.

\begin{figure}
\bigskip
\centertexdraw{
  \drawdim truecm  \linewd.02
  \rhombus \rhombus \rhombus \ldreieck
  \move (-0.866025403784439 -.5)
  \rhombus \rhombus \rhombus \rhombus \ldreieck
  \move (-1.732050807568877 -1)
  \rhombus \rhombus \rhombus \rhombus \rhombus
  \move (-1.732050807568877 -1)
  \rdreieck
  \rhombus \rhombus \rhombus \rhombus
  \move (-1.732050807568877 -2)
  \rdreieck
  \rhombus \rhombus \rhombus
  \move (-0.866025403784439 -1.5)
  \bsegment
    \rlvec(0.866025403784439 -.5) \rlvec(-0.866025403784439 -.5)
	   \rlvec(-0.866025403784439 .5)
    \lfill f:0.8
  \esegment
  \linewd.05
  \move (-1.732050807568877 -2)
  \rlvec(5.196152422706631 0)
\move(8 0)
\bsegment
  \drawdim truecm  \linewd.02
  \rhombus \rhombus \rhombus \ldreieck
  \move (-0.866025403784439 -.5)
  \rhombus \rhombus \rhombus \rhombus \ldreieck
  \move (-1.732050807568877 -1)
  \rhombus \rhombus \rhombus \rhombus \rhombus \ldreieck
  \move (-1.732050807568877 -1)
  \rdreieck \rhombus \rhombus \rhombus \rhombus \rhombus
  \move (-1.732050807568877 -2)
  \rdreieck \rhombus \rhombus \rhombus \rhombus
  \move (-1.732050807568877 -3)
  \rdreieck \rhombus \rhombus \rhombus
  \move (1.732050807568877 -2)
  \bsegment
    \rlvec(0.866025403784439 -.5) \rlvec(-0.866025403784439 -.5)
	   \rlvec(-0.866025403784439 .5)
    \lfill f:0.8
  \esegment
  \linewd.05
  \move (-1.732050807568877 -2.5)
  \rlvec(5.196152422706631 0)
\esegment
\htext (-2 -6){\small A hexagon with sides $N,M,N,$}
\htext (-2 -6.5){\small $N,M,N$ and fixed rhombus }
\htext (-2 -7){\small $l$, where $N=3$, $M=2$, $l=1$.}
\htext (-2 -7.5){\small The thick  horizontal line indicates}
\htext (-2 -8){\small the symmetry axis.}
\htext (5.8 -6){\small A hexagon with sides $N,M,N,$}
\htext (5.8 -6.5){\small $N,M,N$ and fixed rhombus }
\htext (5.8 -7){\small $l$, where $N=3$, $M=3$, $l=2$.}
\htext (5.8 -7.5){\small The thick  horizontal line indicates}
\htext (5.8 -8){\small the symmetry axis.}
\htext (-.9 0.1){$N$}
\htext (2.4 0.1){$N$}
\htext (-2.5 -2){$M$}
\htext (7.1 0.1){$N$}
\htext (10.4 0.1){$N$}
\htext (5.5 -2.6){$M$}
}
\caption{}
\label{fig:1}
\end{figure}

\begin{thm}
\label{thm:MEven}
Let $m$ be a nonnegative integer and $N$ be a positive integer. The number
of rhombus tilings of a hexagon with sides $N,2m,N,N,2m,N$, which
contain the $l$-th rhombus on the symmetry axis 
which cuts through the sides of length $2m$, equals
\begin{multline}
\label{eq:MEven}
\frac{m\binom{m+N}{m}\binom{m+N-1}{m}}{\binom{2m+2N-1}{2m}}
\left(
\sum_{e=0}^{l-1}(-1)^{e}\binom{N}{e}
	\frac{
		(N-2e)\po{\frac{1}{2}}{e}
	}{
		(m+e)(m+N-e)\po{\frac{1}{2}-N}{e}
	}
\right)
\\\times
\prod_{i=1}^N\prod_{j=1}^N\prod_{k=1}^{2m}\frac{i+j+k-1}{i+j+k-2},
\end{multline}
where the shifted factorial
$(a)_k$ is defined by $(a)_k:=a(a+1)\cdots(a+k-1)$,
$k\ge1$, $(a)_0:=1$.
\end{thm}

\begin{thm}
\label{thm:MOdd}
Let $m$ and $N$ be positive integers. The number
of rhombus tilings of a hexagon with sides $N+1,2m-1,N+1,N+1,2m-1,N+1$, which
contain the $l$-th rhombus on the symmetry axis
which cuts through the sides of length $2m-1$, equals
\begin{multline}
\label{eq:MOdd}
\frac{m\binom{m+N}{m}\binom{m+N-1}{m}}{\binom{2m+2N-1}{2m}}
\left(\sum_{e=0}^{l-1}(-1)^{e}\binom{N}{e}
	\frac{
		(N-2e)\po{\frac{1}{2}}{e}
	}{
		(m+e)(m+N-e)\po{\frac{1}{2}-N}{e}
	}
\right)
\\\times
\prod_{i=1}^{N+1}\prod_{j=1}^{N+1}\prod_{k=1}^{2m-1}\frac{i+j+k-1}{i+j+k-2}.
\end{multline}
\end{thm}

In general, the sum in \eqref{eq:MEven} and \eqref{eq:MOdd} (note
that it is indeed exactly the same sum) does not simplify. (The only case
where the sum is known to simplify is for $N=2n-1$, $m=l=n$, see
\cite[Corollary~3]{CiucKratAB} and \cite[Theorem~12, Lemma~13]{HeGeAA}.) Given that this is
the case, a natural ``next" question is to ask what the quantities in
\eqref{eq:MEven} and \eqref{eq:MOdd} are ``roughly", or, more
precisely, what the proportion of rhombus tilings, which contain the
fixed rhombus, in the total number of rhombus tilings is
asymptotically as the hexagon becomes large. Indeed,
from Theorems~\ref{thm:MEven} and \ref{thm:MOdd}, we are able to derive
an ``arcsine law" for this kind of enumeration. 

\begin{thm}\label{thm:Asymptotics}
Let $a$ be any nonnegative real number, let $b$ be a real number with 
$0<b<1$.
For $m\sim a N$ and $l\sim b N$, the proportion of rhombus tilings of a
hexagon with sides $N,2m,N,N,2m,N$ or $N+1,2m-1,N+1,N+1,2m-1,N+1$, 
which contain the
$l$-th rhombus on the symmetry axis which cuts through the sides of length
$2m$, respectively $2m-1$, in the total
number of rhombus tilings is asymptotically
\begin{equation}
\label{eq:Asymptotics}
\frac{2}{\pi}\arcsin\left(
	\frac{\sqrt{b(1-b)}}{\sqrt{(a+b)(a-b+1)}}
\right)
\end{equation}
as $N$ tends to infinity.
\end{thm}

This result is in accordance with the Cohn, Larsen and Propp result
\cite[Theorem~1]{CoLPAA} that was mentioned above. The latter result
does in fact give an asymptotic expression for this kind of
enumeration for an {\em arbitrary} (semiregular) hexagon and an {\em
arbitrary} fixed rhombus. (The result is even much stronger: The
resulting function can even be used as a ``density function"; consult
the paper \cite{CoLPAA} for further information.) Indeed, the value ${\cal
{P}}_{N,2m,N}(0,\sqrt{3}\,(-N/2+l))$, as defined in \cite[Theorem~1]{CoLPAA},
for $m=aN$ and $l=bN$ equals 
$$\frac {1} {\pi}\arccot\left(\frac {a(1+a)-b(1-b)}
{2\sqrt{a(1+a)b(1-b)}}\right),$$
which, in view of the formula
$$2\arcsin x=\arccot\left(\frac {1-2x^2} {2x\sqrt{1-x^2}}\right),$$
is exactly equal to \eqref{eq:Asymptotics}, as it should be.

In the next section we describe proofs of 
Theorems~\ref{thm:MEven}, \ref{thm:MOdd} and \ref{thm:Asymptotics}.
For the proofs of Theorems~\ref{thm:MEven} and \ref{thm:MOdd} we
build on the approach of \cite{CiucKratAB}. The basic ingredients are
an application of Ciucu's Matchings Factorization Theorem 
\cite[Theorem~1.2]{CiucAB}, the standard correspondence between rhombus
tilings and nonintersecting lattice paths, and evaluations of the
determinants which result from the nonintersecting lattice paths.
As opposed to \cite{CiucKratAB}, here we have to
take care of a few subtleties, which do not arise in 
\cite{CiucKratAB} when the fixed rhombus is the central rhombus.
Besides, we offer a few simplifications compared to
\cite{CiucKratAB}.
The proof of a crucial auxiliary lemma is deferred to Section~3.
Finally, in Section~4, we point to further directions in this
research, and list a few open problems and conjectures.

\section{The proofs}
{\sc Proof of Theorems~\ref{thm:MEven} and \ref{thm:MOdd}}.
The proofs of both Theorems are very similar. We will mainly concentrate
on the proof of Theorem~\ref{thm:MEven}.

There are four basic steps. 

{\it Step~1. Application of the Matchings Factorization Theorem}.
First, rhombus tilings of the hexagon with sides $N,2m,N,N,2m,N$
can be interpreted as perfect matchings of the dual graph of the
triangulated hexagon, i.e., the (bipartite) graph $G(V,E)$, where
the set of vertices $V$ consists of the triangles of the hexagon's
triangulation, and where
two vertices are connected by an edge if the corresponding triangles
are adjacent. Enumerating only those rhombus tilings which contain
a fixed rhombus, under this translation amounts to enumerating only
those perfect matchings which contain the edge corresponding to this
rhombus, or, equivalently, we may consider just perfect matchings of the
graph which results from $G(V,E)$ by removing this edge. 
Clearly, since the fixed rhombus was located on the symmetry axis,
this graph is symmetric.
Hence, we may apply
Ciucu's Matchings Factorization
Theorem \cite[Theorem~1.2]{CiucAB}. In general, this theorem says that 
the number of perfect matchings of a symmetric graph $G$ equals a certain
power of 2 times the number of perfect matchings of a graph $G^+$ (which
is, roughly speaking, the ``upper half" of $G$) times a weighted count 
of perfect
matchings of a graph $G^-$ (which is, roughly speaking, the ``lower half" of
$G$), in which the edges on the symmetry axis count with weight $1/2$ only.
Applied to our case, and retranslated into rhombus tilings, the Matchings
Factorization Theorem implies the following:

{\em The number
of rhombus tilings of a hexagon with sides $N,2m,N,N,2m,N$ which
contain the $l$-th rhombus on the symmetry axis 
which cuts through the sides of length $2m$, equals
\begin{equation}\label{eq:MatchFact}
2^{N-1}R(S'(N,m))\tilde R(C(N,m,l)),
\end{equation}}
where $S'(N,m)$ denotes the ``upper half" of our hexagon with the fixed
rhombus removed (see the left half of
Figure~\ref{fig:2}), where $R(S'(m,n))$ denotes 
the number of rhombus tilings of $S'(m,n)$, where $C(N,m,l)$ denotes the 
``lower half" (again, see the left half of
Figure~\ref{fig:2}), and where $\tilde R(C(N,m,l))$
denotes the weighted count of rhombus tilings of $C(N,m,l)$ in which
each of the top-most (horizontal) rhombi counts with weight $1/2$.
(Both, $S'(N,m)$ and $C(N,m,l)$ are roughly pentagonal. The notations $S'(N,m)$
and $C(N,m,l)$ stand for ``simple part" and ``complicated part", respectively,
as it will turn out that the count $R(S'(N,m))$ will be rather 
straight-forward, while the count $\tilde R(C(N,m,l))$ will turn out be
considerably harder.)

It is immediately obvious, that
along the left-most and right-most vertical strip of $S'(N,m)$, the rhombi
are {\em uniquely\/} determined. Hence, we may safely remove these strips
(see the left half of Figure~\ref{fig:2}, the strips are shaded). 
Let us denote
the resulting region by $S(N-1,m)$. From \eqref{eq:MatchFact}
we obtain that
{\em the number
of rhombus tilings of a hexagon with sides $N,2m,N,N,2m,N$, which
contain the $l$-th rhombus on the symmetry axis 
which cuts through the sides of length $2m$, equals
\begin{equation}
\label{eq:CiucuMEven}
2^{N-1}R(S(N-1,m))\tilde R(C(N,m,l)).
\end{equation}}

\begin{figure}
\bigskip
\centertexdraw{
\drawdim truecm  \linewd.02
\move(0 -1)
\bsegment
  \rhombus \rhombus \rhombus \ldreieck
  \move (-0.866025403784439 -.5)
  \rhombus \rhombus \rhombus
  \move (-1.732050807568877 -1)
  \rhombus \rhombus
  \move (-1.732050807568877 -1)
  \rdreieck
  \move (1.732050807568877 -4)
  \rhombus \ldreieck
  \move (-0 -4)
  \rhombus \rhombus \rhombus
  \move (-.866025403784439 -4.5)
  \rhombus \rhombus \rhombus
  \move (-1.732050807568877 -4)
  \rdreieck \rhombus \rhombus \rhombus
  \move (-1.732050807568877 -1)
  \bsegment
    \rlvec(0.866025403784439 .5) \rlvec(0 -1)
	   \rlvec(-0.866025403784439 -.5)
    \lfill f:0.5
  \esegment
  \move (2.598076211353316 -.5)
  \bsegment
    \rlvec(0.866025403784439 -.5) \rlvec(0 -1)
	   \rlvec(-0.866025403784439 .5)
    \lfill f:0.5
  \esegment
  \linewd.08
  \move (-1.732050807568877 -1)
  \RhombusB
  \move (2.598076211353316 -.5)
  \RhombusC
  \htext (-2.0 0.6){
  	$\overbrace{\hphantom{xxxxxxxxxxxxxxxxxxxxxxx}}^{S^\prime(3,1)}$
	}
  \htext (-0.8 -2.7){$\underbrace{\hphantom{xxxxxxxxxxxxxx}}_{S(2,1)}$}
  \htext (-2.0 -7.2){
  	$\underbrace{\hphantom{xxxxxxxxxxxxxxxxxxxxxxx}}_{C(3,1,1)}$
	}
\esegment
\move(7 -1)
\bsegment
  \drawdim truecm  \linewd.02
  \rhombus \rhombus \rhombus \ldreieck
  \move (-0.866025403784439 -.5)
  \rhombus \rhombus \rhombus \rhombus
  \move (-1.732050807568877 -1)
  \rhombus \rhombus \rhombus
  \move (-1.732050807568877 -1)
  \rdreieck \rhombus
  \move (2.598076211353316 -4.5)
  \ldreieck
  \move (1.732050807568877 -5)
  \rhombus \ldreieck
  \move (-0.866025403784439 -4.5)
  \rhombus \rhombus \rhombus \rhombus
  \move (-1.732050807568877 -4)
  \rdreieck \rhombus \rhombus \rhombus \rhombus
  \move (-1.732050807568877 -5)
  \rdreieck \rhombus \rhombus \rhombus
  \move (-1.732050807568877 -4)
  \bsegment
    \rlvec(0.866025403784439 -.5) \rlvec(0 -2)
	   \rlvec(-0.866025403784439 .5)
    \lfill f:0.5
  \esegment
  \move (3.464101615137754 -4)
  \bsegment
    \rlvec(-0.866025403784439 -.5) \rlvec(0 -2)
	   \rlvec(0.866025403784439 .5)
    \lfill f:0.5
  \esegment
  \linewd.08
  \move (-1.732050807568877 -4)
  \RhombusC
  \move (-1.732050807568877 -5)
  \RhombusC
  \move (2.598076211353316 -4.5)
  \RhombusB
  \move (2.598076211353316 -5.5)
  \RhombusB
  \htext (-2.0 0.6){
  	$\overbrace{\hphantom{xxxxxxxxxxxxxxxxxxxxxxx}}^{S(3,1)}$
	}
  \htext (-2.0 -3.8){
  $\overbrace{\hphantom{xxxxxxxxxxxxxxxxxxxxxxx}}^{C^\prime(3,2,2)}$
	}
  \htext (-0.9 -8.2){
  	$\underbrace{\hphantom{xxxxxxxxxxxxxx}}_{C(2,2,2)}$
	}
\esegment
\move (-1.5 -9.7)
\bsegment
	\htext (0 0){\small The hexagons from Figure~\ref{fig:1} cut
along the symmetry axis, according}
	\htext (0 -.5){\small to the Matchings Factorization Theorem. In the
shaded regions,}
	\htext (0 -1){\small the ``forced" (uniquely determined) rhombi are shown with thick lines.}
\esegment
}
\caption{Hexagons, cut in two}
\label{fig:2}
\end{figure}
Similarly, for the case of Theorem~\ref{thm:MOdd}, we obtain that
{\em the number
of rhombus tilings of a hexagon with sides $N,2m-1,N,N,2m-1,N$, which
contain the $l$-th rhombus on the symmetry axis 
which cuts through the sides of length $2m-1$, equals
\begin{equation}
\label{eq:CiucuMOdd}
2^{N-1}R(S(N,m-1))\tilde R(C(N-1,m,l)).
\end{equation}}
(See the right half of Figure~\ref{fig:2}.
Note that in the case of Theorem~\ref{thm:MOdd}, an application of the Matchings Factorization Theorem would directly give us
$2^{N-1}R(S(N,m-1))\tilde R(C'(N,m,l))$, with $C'(N,m,l)$ the region
as indicated in Figure~\ref{fig:2}. However, similarly to before, any
rhombus tiling  of
the ``complicated part'' $C'(N,m,l)$ is {\em uniquely\/} determined in the
left-most and right-most vertical strip of $C'(N,m,l)$. Removing these strips then yields
$C(N-1,m,l)$.)

\smallskip
{\it Step~2. From rhombus tilings to nonintersecting lattice paths.}
There is a standard translation from rhombus tilings to nonintersecting
lattice paths. We apply it to our regions $S(N,m)$ and $C(N,m,l)$. 
Figure~\ref{fig:3} illustrates this translation for the (``complicated")
lower parts in Figure~\ref{fig:2}.

\begin{figure}
\bigskip
\centertexdraw{
  \drawdim truecm  \linewd.02
\move(0 0)
\bsegment
\bsegment
  \move (1.732050807568877 0)
  \rhombus \ldreieck
  \move (-0 0)
  \rhombus \rhombus \rhombus
  \move (-.866025403784439 -0.5)
  \rhombus \rhombus \rhombus
  \move (-1.732050807568877 -0)
  \rdreieck \rhombus \rhombus \rhombus
  \linewd.08
  \move (-1.732050807568877 0)
  \RhombusC \RhombusC
  \rmove(-0.866025403784439 0.5)
  \RhombusA \RhombusB \RhombusA
  \rmove(-0.866025403784439 2.5)
  \RhombusA \RhombusA \RhombusB
  \rmove (0 2)
  \RhombusA \RhombusB
  \linewd.05
  \move (-0.4330127018922194 -.25)
  \lcir r:.1
  \hdSchritt \vdSchritt \hdSchritt
  \lcir r:.1
  \move (0.4330127018922194 .25)
  \lcir r:.1
  \hdSchritt \hdSchritt \vdSchritt
  \lcir r:.1
  \move (2.165063509461096 .25)
  \lcir r:.1
  \hdSchritt \vdSchritt
  \lcir r:.1
\esegment
\move(7 0)
\bsegment
  \linewd.02
  \move (1.732050807568877 -0.5)
  \ldreieck
  \move (-0.866025403784439 -0)
  \rhombus \rhombus \rhombus \ldreieck
  \move (-0.866025403784439 -0)
  \rdreieck \rhombus \rhombus \rhombus
  \move (-0.866025403784439 -1)
  \rdreieck \rhombus \rhombus
  \linewd.08
  \move (-0.866025403784439 0)
  \RhombusC
  \rmove (-0.866025403784439 -0.5)
  \RhombusC
  \move (-0.866025403784439 0)
  \RhombusA \RhombusB \RhombusB \RhombusA
  \move (0.866025403784439 0)
  \RhombusC
  \rmove (-0.866025403784439 -0.5)
  \RhombusC
  \move (1.732050807568877 -0.5)
  \RhombusB \RhombusB
  \linewd.05
  \move (-0.4330127018922194 .25)
  \lcir r:.1
  \hdSchritt \vdSchritt \vdSchritt \hdSchritt
  \lcir r:.1
  \move (2.165063509461096 -.25)
  \lcir r:.1
  \vdSchritt \vdSchritt
  \lcir r:.1
\esegment
\esegment
\move(0 -6.7)
\bsegment
	\move(-1 0)
	\bsegment
		\linewd.02
		\move (0 0) \lvec (4.5 0)
		\move (0 0) \lvec (0 3.5)
		\linewd.01
		\move (0 1) \rlvec (4.5 0)
		\move (0 2) \rlvec (4.5 0)
		\move (0 3) \rlvec (4.5 0)
		\move (1 0) \rlvec (0 3.5)
		\move (2 0) \rlvec (0 3.5)
		\move (3 0) \rlvec (0 3.5)
		\move (4 0) \rlvec (0 3.5)
		\linewd.07
		\move (0 1) \lcir r:.1 \rlvec(1 0) \rlvec(0 -1) \rlvec(1 0)  
\lcir r:.1
		\move (1 2) \lcir r:.1 \rlvec(2 0) \rlvec(0 -1) \lcir r:.1
		\move (3 3) \lcir r:.1 \rlvec(1 0) \rlvec(0 -1) \lcir r:.1
	\esegment
	\move(5 0)
	\bsegment
		\linewd.02
		\move (0 0) \lvec (4.5 0)
		\move (0 0) \lvec (0 3.5)
		\linewd.01
		\move (0 1) \rlvec (4.5 0)
		\move (0 2) \rlvec (4.5 0)
		\move (0 3) \rlvec (4.5 0)
		\move (1 0) \rlvec (0 3.5)
		\move (2 0) \rlvec (0 3.5)
		\move (3 0) \rlvec (0 3.5)
		\move (4 0) \rlvec (0 3.5)
		\linewd.07
		\move (0 2) \lcir r:.1 \rlvec(1 0) \rlvec(0 -2) \rlvec(1 0) \lcir r:.1
		\move (3 3) \lcir r:.1 \rlvec(0 -2) \lcir r:.1
	\esegment
\esegment
\move(-3.5 -8.5)
\bsegment
	\htext (0 0){\small Tilings for the ``complicated parts'' from
		Figure~\ref{fig:2}, interpreted as lattice paths.}
\esegment
}
\caption{Lattice path interpretation}
\label{fig:3}
\end{figure}

For the ``simple" pentagonal part $S(N,m)$ we obtain the following:
The number\break $R(S(N,m))$ 
of rhombus tilings of $S(N,m)$ equals the number of families
$(P_1,P_2,\dots,P_N)$ of nonintersecting lattice paths consisting of
horizontal unit steps in the positive direction and vertical unit steps
in the negative direction, where $P_i$ runs from $(i,i)$ to 
$(2i, i-m)$,
$i=1,2,\dots,N$.

Similarly, for the ``complicated" pentagonal part $C(N,m,l)$ we obtain:
The weighted count $\tilde R(C(N,m,l))$ 
of rhombus tilings of $C(N,m,l)$ equals the weighted count of families
$(P_1,P_2,\dots,P_N)$ of nonintersecting lattice paths consisting of
horizontal unit steps in the positive direction and vertical unit steps
in the negative direction, where $P_i$ runs from 
$(2i-N-1,i+m)$ to $(i,i)$ if $i\neq l$, while
$P_l$ runs from $(2l-N,l+m)$ to $(l,l)$, with the additional twist
that for $i\neq l$ path $P_i$ has weight $1/2$ if it starts with a
horizontal step.

\smallskip
{\it Step~3. From nonintersecting lattice paths to determinants}.
Now, by using the main theorem on nonintersecting lattice paths 
\cite[Corollary~2]{GeViAB} (see also \cite[Theorem~1.2]{StemAE}),
we may write $R(S(N,m))$ and $\tilde R(C(N,m,l))$ as
determinants. Namely, 
we have
\begin{align}
R(S(N,m))	
&=\det_{1\leq i,j\leq N}\left(\binom{m+i}{m-i+j}\right)\label{eq:SNm1},
\end{align}
and
\begin{equation}
\label{eq:CNml1}
\tilde R(C(N,m,l))	=
	\det_{1\leq i,j\leq N}\left(
		\begin{cases}
\frac{(N+m-i)!}{(m+i-j)!\,(N+j-2i+1)!}(m+\frac{N-j+1}{2})&\text{ if }i\neq l\\
			\frac{(N+m-i)!}{(m+i-j)!\,(N+j-2i)!}&\text{ if }i=l
		\end{cases}
	\right).
\end{equation}

\smallskip
{\it Step~4. Determinant evaluations}. Clearly, once we are able
to evaluate the determinants in \eqref{eq:SNm1} and \eqref{eq:CNml1},
Theorems~\ref{thm:MEven} and \ref{thm:MOdd} will immediately follow
from \eqref{eq:CiucuMEven} and \eqref{eq:CiucuMOdd}, respectively,
upon routine simplification.
Indeed, for the determinant in \eqref{eq:SNm1} we have the following.
\begin{lem}
\label{lem:simplepart}
\begin{align}
&\det_{1\leq i,j\leq N}\left(\binom{m+i}{m-i+j}\right)
 =\prod_{i=1}^{N}\frac{(N+m-i+1)!\,(i-1)!\,\po{2m+i+1}{i-1}}{(m+i-1)!\,(2N-2i+1)!}.
			\label{eq:SNm2}
\end{align}
\end{lem}
\begin{pf}This determinant evaluation follows without difficulty from
a determinant lemma in \cite[Lemma~2.2]{KratAM}. The corresponding
computation is contained in the proof of Lemma~9 in \cite{CiucKratAB}
(our determinant becomes the same as in \cite{CiucKratAB} when the
order of rows and columns is reversed),
and also in the proof of Theorem~5 in \cite{KratAK}
(set $r=N$, $\lambda_s=m$, $B=2$, $a+\alpha-b=2m$ there, and then reverse
the order of rows and columns). 
\end{pf}

On the other hand, the determinant in \eqref{eq:CNml1} evaluates as follows.
\begin{lem}
\label{lem:complexpart}
\begin{multline}
\label{eq:CNml2}
	\det_{1\leq i,j\leq N}\left(
		\begin{cases}
\frac{(N+m-i)!}{(m+i-j)!\,(N+j-2i+1)!}(m+\frac{N-j+1}{2})&\text{ if }i\neq l\\
			\frac{(N+m-i)!}{(m+i-j)!\,(N+j-2i)!}&\text{ if }i=l
		\end{cases}
	\right)\\
= \prod_{i=1}^{N}\frac{(N+m-i)!}{(m+i-1)!\,(2N-2i+1)!}
\prod _{i=1}^{\floor{N/2}}
	\left(
		\po{m+i}{N-2i+1}\,\po{m+i+\frac{1}{2}}{N-2i}
	\right)\kern1.5cm
\\
\times
2^{\frac{(N-1)(N-2)}{2}}\frac{
	\po{m}{N+1}\prod_{j=1}^{N}(2j-1)!
}{
	N!\prod_{i=1}^{\floor{\frac{N}{2}}}\po{2i}{2N-4i+1}
}
\sum_{e=0}^{l-1}
(-1)^{e}\binom{N}{e}\frac{
	(N-2e)\,\po{\frac{1}{2}}{e}
}{
	(m+e)\,(m+N-e)\,\po{\frac{1}{2}-N}{e}
}.
\end{multline}
\end{lem}
This determinant evaluation is much more complex than the determinant
evaluation of Lemma~\ref{lem:simplepart}, and, as such, is the most
difficult part in our derivation of Theorems~\ref{thm:MEven} and
\ref{thm:MOdd}. We defer the proof of Lemma~\ref{lem:simplepart}
to the next section.

\smallskip
Altogether, Steps~1--4 establish Theorems~\ref{thm:MEven} and
\ref{thm:MOdd}.
\hfil\qed

\bigskip
{\sc Proof of Theorem~\ref{thm:Asymptotics}}.
{}From MacMahon's formula~\eqref{eq:MacMahon} for the total number of rhombus
tilings together with Theorems~\ref{thm:MEven} and \ref{thm:MOdd} we deduce
immediately that the proportion is indeed the same for both cases
$N,2m,N,N,2m,N$ and $N+1,2m-1,N+1,N+1,2m-1,N+1$, and that it is given by
\begin{equation}
\label{eq:proportion}
\frac{m\binom{m+N}{m}\binom{m+N-1}{m}}{\binom{2m+2N-1}{2m}}
\sum_{e=0}^{l-1}(-1)^{e}\binom{N}{e}
	\frac{
		(N-2e)\,\po{\frac{1}{2}}{e}
	}{
		(m+e)\,(m+N-e)\,\po{\frac{1}{2}-N}{e}
	}.
\end{equation}

We write the sum in \eqref{eq:proportion} in a hypergeometric fashion,
to get
\begin{equation}
\label{eq:proportion-hyp-0}
\frac{(2N-1)!\,\left(\po{m+1}{N-1}\right)^2}{(N-1)!^2\,\po{2m+1}{2N-1}}
\sum _{e=0} ^{l-1}\frac {(-N)_e\,(1-\frac{N}{2})_e\,
(m)_e\,(-m-N)_e\,(\frac{1}{2})_e} 
{(-\frac{N}{2})_e\,(1-m-N)_e\,(1+m)_e\,(\frac{1}{2}-N)_e\,e!}.
\end{equation}
Next we shall apply Whipple's transformation (see \cite[(2.4.1.1)]{SlatAC}),
which reads
\begin{multline}
\label{eq:Whipple}
\HypF{7}{6}{a,1+\frac{a}{2},b,c,d,e,-n}%
{\frac{a}{2},1+a-b,1+a-c,1+a-d,1+a-e,1+a+n}{1}\\
=\frac{\po{a+1}{n}\po{a-d-e+1}{n}}{\po{a-d+1}{n}\po{a-e+1}{n}}
\HypF{4}{3}{a-b-c+1,d,e,-n}{a-b+1,a-c+1,-a+d+e-n}{1}.
\end{multline}
Here, we used the standard hypergeometric notation
\begin{equation}
\HypF{r}{s}{a_1,\dots,a_r}{b_1,\dots,b_s}{z} =
	\sum_{k=0}^\infty
		\frac{
			\po{a_1}{k}\cdots\po{a_r}{k}
		}{
			k!\po{b_1}{k}\cdots\po{b_s}{k}
		}z^k.
\end{equation}
By setting $a=-N,b=m,c=-m-N,d=\frac{1}{2},e=-N+l,n=l-1$ in \eqref{eq:Whipple},
we obtain as a limit case the following transformation:
\begin{multline}
\sum _{e=0} ^{l-1}\frac {(-N)_e\,(1-\frac{N}{2})_e\,
(m)_e\,(-m-N)_e\,(\frac{1}{2})_e} 
{(-\frac{N}{2})_e\,(1-m-N)_e\,(1+m)_e\,(\frac{1}{2}-N)_e\,e!}\\
=\frac{\po{-N+1}{l-1}\po{-l+\frac{1}{2}}{l-1}}%
{\po{-N+\frac{1}{2}}{l-1}\po{-l+1}{l-1}}
\HypF{4}{3}{1,\frac{1}{2},l-N,1-l}{1+m,1-m-N,\frac{3}{2}}{1}.
\end{multline}
Thus, expression \eqref{eq:proportion-hyp-0} turns into
\begin{equation}\notag
\frac{(2N-1)!\,\left(\po{m+1}{N-1}\right)^2}{(N-1)!^2\,\po{2m+1}{2N-1}}
\frac{\po{-N+1}{l-1}\,\po{-l+\frac{1}{2}}{l-1}}%
{\po{-N+\frac{1}{2}}{l-1}\,\po{-l+1}{l-1}}
\HypF{4}{3}{1,\frac{1}{2},l-N,1-l}{1+m,1-m-N,\frac{3}{2}}{1}.
\end{equation}

Next we apply Bailey's transformation (see \cite[(4.3.5.1)]{SlatAC}) between
two balanced $\HypFsimple{4}{3}$-series,
\begin{multline}
\HypF{4}{3}{a,b,c,-n}%
{e,f,1+a+b+c-e-f-n}{1}=
\frac{\po{e-a}{n}\po{f-a}{n}}%
{\po{e}{n}\po{f}{n}}\times\\
\HypF{4}{3}{-n,a,a+c-e-f-n+1,a+b-e-f-n+1}%
{a+b+c-e-f-n+1,a-e-n+1,a-f-n+1}{1}.
\end{multline}
with $a=1,b=1/2,c=l-N,e=m+1,f=-m-N+1,n=l-1$.
This gives
\begin{multline}\label{eq:fact4F3}
\frac{(2l)!\,(2m)!\,(m+N-1)!\,(m+N)!\,(2N-2l+2)!}%
{4(l+m-1)(m+N-l+1)(l-1)!\,l!\,(m-1)!}\\
\times
\frac{1}{m!\,(N-l)!\,(N-l+1)!\,(2m+2N-1)!}
\HypF{4}{3}{1-l,1,1,\frac{3}{2}-l+N}%
{\frac{3}{2},2-l-m,2-l+m+N}{1}
\end{multline}
for the ratio \eqref{eq:proportion}.
Now we substitute $m\sim a N$ and $l\sim b N$ and perform the limit
$N\rightarrow\infty$. With Stirling's formula we determine the
limit for the quotient of factorials in front of the $_4F_3$-series
in \eqref{eq:fact4F3} as
${2\sqrt{a(a+1)}\sqrt{b(1-b)}}/({\pi(a-b+1)(a+b)})$.
For the $\HypFsimple{4}{3}$-series itself, 
we may exchange limit and summation
by uniform convergence:
\begin{equation}\notag
\lim_{N\rightarrow\infty}\HypF{4}{3}{1-l,1,1,\frac{3}{2}-l+N}%
{\frac{3}{2},2-l-m,2-l+m+N}{1}=
\HypF{2}{1}{1,1}{\frac{3}{2}}{ {(1-b)b\over(a-b+1)(a+b)} }.
\end{equation}
A combination of these results and use of the identity
(see \cite[p.~463, (133)]{Prudnikov})
\begin{equation}\notag
\HypF{2}{1}{1,1}{\frac{3}{2}}{z}=\frac{\arcsin\sqrt{z}}{\sqrt{z(1-z)}}
\end{equation}
finish the proof.
\hfill\qed

\section{Proof of Lemma~\ref{lem:complexpart}}
The method that we use for this proof is also applied 
successfully in \cite{KratBD,CiucKratAB,KratBG,KratBH,KratBI,KrZeAA}
(see in particular the tutorial description in \cite[Sec.~2]{KratBI}).

First of all, we take appropriate factors out of the determinant
in \eqref{eq:CNml2}.
To be precise, we take 
$$
\frac{(N+m-i)!}{(m+i-1)!\,(2N-2i+1)!}
$$
out of the $i$-th row of the determinant, $i=1,2,\dots,N$. Thus we obtain
\begin{multline} \label{eq:detdef}
 \prod_{i=1}^{N}\frac{(N+m-i)!}{(m+i-1)!\,(2N-2i+1)!}\\\times
	\det_{1\leq i,j\leq N}\left(
		\begin{cases}
			\po{m+i-j+1}{j-1}\po{N+j-2i+2}{N-j}\frac{N+2m-j+1}{2}
				&\text{ if } i\neq l \\
   		\po{m+i-j+1}{j-1}\po{N+j-2i+1}{N-j+1}
				&\text{ if } i=l
		\end{cases}
	\right)
\end{multline}
for the determinant in \eqref{eq:CNml2}.
Let us denote by $D(m;N,l)$ the $N\times N$-matrix underlying the
determinant in \eqref{eq:detdef}, i.e., the $(i,j)$-entry of $D(m;N,l)$
is given by
\begin{equation}
\label{eq:def-D}
D(m;N,l)_{i,j}:=
	\begin{cases}
			\po{m+i-j+1}{j-1}\po{N+j-2i+2}{N-j}\frac{N+2m-j+1}{2}
				&\text{ if } i\neq l, \\
   		\po{m+i-j+1}{j-1}\po{N+j-2i+1}{N-j+1}
				&\text{ if } i=l.
	\end{cases}
\end{equation}

Comparison of \eqref{eq:CNml2} and 
\eqref{eq:detdef} yields that \eqref{eq:CNml2} will be proved once
we are able to establish the determinant evaluation
\begin{multline}
\det\left(D(m;N,l)\right)=
\prod _{i=1}^{\floor{N/2}}
	\left(
		\po{m+i}{N-2i+1}\po{m+i+\tfrac{1}{2}}{N-2i}
	\right)
\\\times
2^{\frac{(N-1)(N-2)}{2}}\frac{
	\po{m}{N+1}\prod_{j=1}^{N}(2j-1)!
}{
	N!\prod_{i=1}^{\floor{\frac{N}{2}}}\po{2i}{2N-4i+1}
} 
\sum_{e=0}^{l-1}
(-1)^{e}\binom{N}{e}\frac{
	(N-2e)\po{\frac{1}{2}}{e}
}{
	(m+e)(m+N-e)\po{\frac{1}{2}-N}{e}
}.\label{eq:main}
\end{multline}

\medskip

For the proof of \eqref{eq:main} we proceed in several steps (see below).
An outline is as
follows. In the first step we show that 
$\prod _{i=1}^{\floor{N/2}}\po{m+i}{N-2i+1}$ is a factor of
$\det\left(D(m;N,l)\right)$  
as a polynomial in $m$. In the second step we
show that 
$\prod_{i=1}^{\floor{N/2}}\po{m+i+\frac{1}{2}}{N-2i}$
is a factor of $\det\left(D(m;N,l)\right)$.
In the third step we determine the maximal degree of
$\det\left(D(m;N,l)\right)$ as a
polynomial in $m$, which turns out to be $\binom {N+1}2-1$. From a
combination of these three steps we are forced to conclude that
\begin{equation} \label{eq:polydef}
\det\left(D(m;N,l)\right)=\prod _{i=1}^{\floor{N/2}}
	\left(
		\po{m+i}{N-2i+1}\po{m+i+\tfrac{1}{2}}{N-2i}
	\right)
P(m;N,l),
\end{equation}
where $P(m;N,l)$ is a polynomial in $m$ of degree at most $N-1$.
Finally, in the fourth step, we evaluate $P(m;N,l)$ at $m=0,-1,\dots,-N$.
Namely, for $m=0,-1,\dots,-\floor{N/2}$ we show that
\begin{multline} \label{eq:polyeval}
P(m;N,l)=
(-1)^{m N+(m^2-m)/2}2^{(m^2+m)/2-N+1}\po{m}{m}
\\ \times
\frac{
		\prod_{j=1}^{N-m}(2j-1)!
\prod_{k=1}^{m}(k-1)!^2(N+k-2m-1)!\po{\frac{m-k+1}{2}}{k-1}\po{k-N}{N-m}
}{
		\prod_{i=1}^{m}(N-m-i)!(m-i)!
		\prod_{i=m+1}^{\floor{N/2}}\po{i-m}{N-2i+1}
		\prod_{i=1}^{\floor{N/2}}\po{i-m+\frac{1}{2}}{N-2i}
}.
\end{multline}
Moreover, we show that $P(m;N,l)=P(-N-m;N,l)$, which in combination with 
\eqref{eq:polyeval} gives the evaluation of $P(m;N,l)$ at $m=-\floor{N/2}-1,
\dots,-N+1,-N$.
Clearly, this determines a polynomial of maximal
degree $N-1$ uniquely. In fact, an explicit expression for $P(m;N,l)$ can 
immediately be written down using Lagrange interpolation. 
As it turns out, the resulting expression for $P(m;N,l)$ is exactly
the second line of \eqref{eq:main}.
In view of \eqref{eq:polydef},
this would establish \eqref{eq:main} and, hence, 
finish the proof of the Lemma.

\smallskip
Before going into details of these steps, however, it is useful to record two
auxiliary facts. The reader may, at this point, directly jump to
Steps~1--4,
and come back to the auxiliary facts when they are needed there.

\medskip
{\it Auxiliary Fact I.} There holds the symmetry
\begin{equation}
\label{eq:simple-symmetry}
\det\left(D(m;N,l)\right)=\det\left(D(m;N,N+1-l)\right).
\end{equation}
This symmetry follows immediately from the combinatorial ``origin" of 
the determinant. For, trivially, the number of rhombus tilings which
contain the $l$-th rhombus on the symmetry axis is the same as the number
of rhombus tilings which
contain the $(N+1-l)$-th rhombus. The manipulations that finally lead to
the determinant $\det\left(D(m;N,l)\right)$ do not affect this symmetry,
therefore  $\det\left(D(m;N,l)\right)$ inherits the symmetry.

This symmetry is very useful for our considerations, because for any
claim that we want to prove (and which also obeys this symmetry) we may
freely assume $1\leq l\leq \floor{\frac{N+1}{2}}$ or
$\floor{\frac{N+1}{2}}\leq l\leq N$, whatever is more convenient.

\medskip
{\it Auxiliary Fact II.} There holds the symmetry
\begin{equation}\label{eq:symmetry}
\det\left(D(-N-m;N,l)\right)=(-1)^{\binom {N+1}2-1}
\det\left(D(m;N,l)\right).
\end{equation}
In order to establish \eqref{eq:symmetry}, we claim that
$D(m;N,l)\cdot\left((-1)^j\binom{j-1}{i-1}\right)_{1\le i,j\le N}$ equals $D(-N-m;N,l)$,
except that all the entries in row $l$ have opposite sign. Let us write 
$R(N)$ for the matrix $\left((-1)^j\binom{j-1}{i-1}\right)_{1\le i,j\le N}$.
Since $\det\left(R(N)\right)=(-1)^{\binom{N+1}2}$, this would establish
\eqref{eq:symmetry}.

In order to establish this claim, we have to compute the $(i,j)$-entry in
$D(m;N,l)\cdot R(N)$. For $i\neq l$, we have to show
\begin{multline}
\label{eq:A-symmetry}
(-1)^j\sum_{k=1}^{j}\binom{j-1}{k-1}
	\po{m+i-k+1}{k-1}\po{N+k-2i+2}{N-k}\tfrac{(N+2m-k+1)}{2}\\
=	\po{-N-m+i-j+1}{j-1}\po{N+j-2i+2}{N-j}\tfrac{(-2m-N-j+1)}{2},
\end{multline}
and for $i=l$ we have to show
\begin{multline}
\label{eq:B-symmetry}
(-1)^j\sum_{k=1}^{j}\binom{j-1}{k-1}
	\po{m+i-k+1}{k-1}\po{N+k-2i+1}{N-k+1}\\
=	-\po{-N-m+i-j+1}{j-1}\po{N-2i+j+1}{N-j+1}.
\end{multline}

Note that for $j=1$, equation \eqref{eq:A-symmetry} collapses to
\begin{equation*}
-\po{N-2i+3}{N-1}\tfrac{(2m+N)}{2}=\po{N-2i+3}{N-1}\tfrac{(-2m-N)}{2},
\end{equation*}
which is of course true, so we may assume $j>1$ in the following.

We convert the left-hand side of \eqref{eq:A-symmetry} into hypergeometric
form, to obtain
\begin{equation}
\label{eq:A-symmetry-hyp}
{(-1)^{j-1}\tfrac {1} {2}(-2m-N)\po{N-2i+3}{N-1}}
	\HypF{3}{2}{-2m-N+1,1-j,1-i-m}{-2m-N,3-2i+N}{1}.
\end{equation}
Next we apply the contiguous relation
\begin{equation}
\label{eq:contiguous-relation}
\HypF{3}{2}{a,A_1,A_2}{B_1,B_2}{z} = 
		\HypF{3}{2}{a-1,A_1,A_2}{B_1,B_2}{z}+
   z\frac{A_1 A_2}{B_1B_2}\HypF{3}{2}{a,A_1+1,A_2+1}{B_1+1,B_2+1}{z}
\end{equation}
to the $\HypFsimple{3}{2}$-series in \eqref{eq:A-symmetry-hyp}. 
We want to apply the case where $a=-2m-N+1$ and $B_1=-2m-N$. By inspection, 
in this case parameters cancel inside the two $\HypFsimple{3}{2}$-series on
the right-hand side of the contiguous relation, thus leaving two
$\HypFsimple{2}{1}$-series. So we obtain
\begin{equation}
\HypF{2}{1}{1-i-m,1-j}{3-2i+N}{1} + 
	\frac{(j-1)(1-i-m)}{(-2m-N)(N-2i+3)}\HypF{2}{1}{2-i-m,2-j}{4-2i+N}{1}.
\end{equation}
Each of the two $\HypFsimple{2}{1}$-series can be evaluated by means of the
Chu--Vandermonde-summation (see \cite[(1.7.7), Appendix (III.4)]{SlatAC}),
\begin{equation}
\label{eqEN}
\HypF{2}{1}{a,-n}{ c}{1} = \frac{\po{c-a}{n}}{\po{c}{n}},
\end{equation}
where $n$ is a nonnegative integer. We have to apply the case where $n=j-1$
and $n=j-2$, respectively. Since $j>1$ in our case, the
nonnegativity-condition is
satisfied, and we obtain
\begin{multline}
\HypF{3}{2}{1-2m-N,1-j,1-i-m}{-2m-N,3-2i+N}{1}\\
=	\frac{
		(-N+i-m-1)(N+j+2m-1)\po{N-i+m+2}{j-2}
	}{
		(-2m-N)\po{N-2i+3}{j-1}
	}.
\end{multline}
Inserting this into \eqref{eq:A-symmetry-hyp} shows that \eqref{eq:A-symmetry}
is true.

In order to show \eqref{eq:B-symmetry}, we convert the left-hand side 
into hypergeometric form, to obtain
\begin{equation}
\label{eq:B-symmetry-hyp}
(-1)^{j}\po{N-2i+2}{N}\HypF{2}{1}{1-i-m,1-j}{2-2i+N}{1}.
\end{equation}
Again, we can evaluate this $\HypFsimple{2}{1}$-series by means of 
Chu--Vandermonde summation \eqref{eqEN}. This proves \eqref{eq:B-symmetry}.

\medskip

Now we are ready for heading into the details of Steps~1--4.

\bigskip
{\it Step 1. $\prod _{i=1}^{\floor{N/2}}\po{m+i}{N-2i+1}$ is a 
factor of $\det\left(D(m;N,l)\right)$}. Here, for the first time,
we make use of the
symmetry \eqref{eq:simple-symmetry}. It implies, that we may restrict 
ourselves to $1\leq l\leq \floor{\frac{N+1}{2}}$.

For $i$ between $1$ and $\floor{N/2}$ 
let us consider row $N-i+1$ of the matrix $D(m;N,l)$. We see that 
the $j$-th entry in this row has the form
$$(m+N-i-j+2)_{j-1}\,(-N+2i+j)_{N-j}\tfrac {N+2m-j+1} {2}.$$
Since $(-N+2i+j)_{N-j}=0$ for $j=1,2,\dots,N-2i$, the first
$N-2i$ entries in this row vanish. Therefore $(m+i)_{N-2i+1}$ is a
factor of each entry in row $N-i+1$, $i=1,2,\dots,\floor{N/2}$. Hence, the
complete product $\prod _{i=1}^{\floor{N/2}}\po{m+i}{N-2i+1}$ divides $\det\left(D(m;N,l)\right)$.

\smallskip
{\it Step 2. $\prod_{i=1}^{\floor{N/2}}\po{m+i+\frac{1}{2}}{N-2i}$
is a factor of $D(m;N,l)$}. Again we make use of the symmetry 
\eqref{eq:simple-symmetry}, which allows us to restrict 
ourselves to $1\leq l\leq \floor{\frac{N+1}{2}}$. 

We observe that the product can be rewritten as
$$\prod_{i=1}^{\floor{N/2}}\po{m+i+\tfrac{1}{2}}{N-2i}=
\prod _{e=1} ^{N-2}(m+e+\tfrac {1} {2})^{\min\{e,N-e-1\}}.$$
Therefore, 
because of the other symmetry \eqref{eq:symmetry}, it suffices to
prove that $(m+e+1/2)^{e}$ divides $D(m;N,l)$ for $e=1,2,\dots,\floor{N/2}-1$.
In order to do so, we claim that for each such $e$ there are $e$ 
linear combinations of the columns, which are themselves linearly
independent, that vanish for $m=-e-1/2$. More precisely, we claim that
for $k=1,2,\dots,e$ there holds
\begin{multline} \label{eq:lincomb}
\sum _{j=1} ^{k}\binom kj \cdot(\text {column
$(N-2e+k+j)$ of $D(-e-1/2;N,l)$})\\
-\frac{\po{N-e-l+\frac{1}{2}}{k}}{(-4)^{k}\po{N-e-l+1}{k}}
\cdot(\text {column
$(N-2e)$ of $D(-e-1/2;N,l)$})=0.
\end{multline}
(Note that this operation is really feasible for these values of $e$ and $k$.)
As is not very difficult to see (cf\@. \cite[Sec.~2]{KratBI}) this would
imply that $(m+e+1/2)^e$ divides $D(m;N,l)$.

Obviously, in order to prove \eqref{eq:lincomb} we have to show
\begin{multline}
\label{eq:Atrans}
\sum_{j=0}^{k}\binom{k}{j}
	\po{e+i-j-k-\tfrac{1}{2}-N}{N-2e+j+k}\times\\
	\po{2N+3-2e-2i+j+k}{2e-j-k-1}
	\tfrac{(-j-k-1)}{2}=0
\end{multline}
which is \eqref{eq:lincomb} restricted 
to the $i$-th row, $i\ne l$, and 
\begin{multline}
\label{eq:Btrans}
\sum_{j=0}^{k}\binom{k}{j}
	\po{e+1+l-j-k-\tfrac{1}{2}-N}{N-2e-1+j+k}
	\po{2N-2e-2l+j+k+1}{2e+1-j-k}-\\
	\frac{\po{N-e-l-\frac{1}{2}}{k}}{(-4)^{k}\po{N-e-l+1}{k}}\\
\times\po{e+1+l-j-\tfrac{1}{2}-N}{N-2e-1+j}
	\po{2N-2e-2l+j+1}{2e+1-j} = 0
\end{multline}
which is \eqref{eq:lincomb} restricted to the $l$-th row.

Both equations \eqref{eq:Atrans} and
\eqref{eq:Btrans} can be shown by the same kind of ``hypergeometrics''
(contiguous relation and Chu-Vandermonde) as was used for 
establishing \eqref{eq:symmetry}.

\smallskip
{\it Step 3. $\det\left(D(m;N,l)\right)$ is a polynomial in $m$ of maximal degree
$\binom {N+1}2-1$}.
Clearly, the degree in $m$ of the $(i,j)$-entry in
the determinant $D(m;N,l)$ is $j$ for $i\ne l$, while it is $j-1$ for
$i=l$. Hence, in the defining expansion of the determinant, each term
has degree $\left(\sum _{j=1} ^{N}j\right)-1=\binom {N+1}2-1$.

\smallskip
{\it Step 4. Evaluation of $P(m;N,l)$ at $m=0,-1,\dots,-N$}.
Again, we make use of the symmetry \eqref{eq:simple-symmetry}, 
and this time restrict ourselves to $\floor{\frac{N+1}{2}}\leq l\leq N$.
On the other hand, by the symmetry \eqref{eq:symmetry} and by
the definition \eqref{eq:polydef} of $P(m;N,l)$, we have
$P(m;N,l)=P(-N-m;N,l)$. Therefore, it suffices to
compute the evaluation of $P(m;N,l)$ at $m=0,-1,\dots,-\floor{N/2}$.

What we would like to do is, for any $e$ with $0\le e\le \floor{N/2}$,
to set $m=-e$ in
\eqref{eq:polydef}, compute $\det\left(D(-e;N,l)\right)$,
and then express $P(-e;N,l)$ as
the ratio of $\det\left(D(-e;N,l)\right)$ and the right-hand side product evaluated
at $m=-e$. Unfortunately, this is typically a ratio $0/0$ and, hence,
undetermined. So, we have to first divide both sides of \eqref{eq:polydef}
by the appropriate power of $(m+e)$, and only then set $m=-e$.

Let $e$, $0\le e\le \floor{N/2}$, be fixed.
For $k=0,1,\dots,{e-1}$ we add 
\begin{equation}
\label{eq:colOp}
\sum_{i=1}^k\binom{k}{i}\cdot(\text {column
$(N+1-2e+k+i)$ of $D(m;N,l)$})
\end{equation}
to column ${N+1-2e+k}$ of $D(m;N,l)$. 
The effect is that then $(m+e)$ is a factor
of each entry in column ${N+1-2e+k}$. This can be proven by exactly the same
``hypergeometrics'' as we used for the proofs of Auxiliary Fact~II
and Step~2.  So, we take $(m+e)$ out of
each entry of column ${N+1-2e+k}$, $k=0,1,\dots,{e-1}$ and denote the
resulting matrix by $D_1(m;N,l,e)$. We obtain
\begin{equation}
\label{eq:def-D2-1}
D_1(m;N,l,e)_{i,j}:=
	\begin{cases}
			\po{m+i-j+1}{j-1}\po{N+j-2i+2}{N-j}\frac{N+2m-j+1}{2}
				&\text{ if } i\neq l\\
   		\po{m+i-j+1}{j-1}\po{N+j-2i+1}{N-j+1}
				&\text{ if } i=l
	\end{cases}
\end{equation}
if $j\leq N-2e$ or $j>N-e$, and
\begin{multline}
\label{eq:def-D2-2}
D_1(m;N,l,e)_{i,j}\\
:=	\begin{cases}
			\po{2e+i-k+m-N}{N-2e+k}\po{N-i+m+1}{k}\times& \\
				\quad\quad\po{2N-2e-2i+2k+3}{2e-2k-1}
				&\text{ if } i\neq l\\
			\po{2N-2e-2i+2k+2}{2e-2k}\po{2e+i-k+m-N}{N-e-i+k}\times& \\
				\quad\quad\po{m+e+1}{i-e-1}\po{1-i+m+N}{k}
				&\text{ if } i=l
	\end{cases}
\end{multline}
if $N-2e+1\leq j\leq N-e$.

From what we did so far,
it is straight-forward that we must have
$$D(m;N,l)=(m+e)^e D_1(m;N,l,e).$$
A combination with \eqref{eq:polydef} gives that
\begin{multline} \label{eq:polydet}
P(m;N,l)\\
=\det\left(D_1(m;N,l,e)\right)\prod _{i=1}^{\floor{N/2}}
	\left(
		\po{m+i}{e-i}\po{m+e+1}{N-i-e}\po{m+i+\tfrac{1}{2}}{N-2i}
	\right)^{-1}.
\end{multline}
Now, in this equation, we are able to set $m=-e$.
Hence, in order to determine the evaluation of $P(m;N,l)$ at $m=-e$,
we need the evaluation of $\det\left(D_1(m;N,l,e)\right)$ at $m=-e$. 

Assuming $l\geq\frac{N+1}{2}$, we claim that the following is true:
If $e\geq N+1-l$, then we have
\begin{equation}
\label{eq:detD2zero}
\det\left(D_1(-e;N,l,e)\right)=0,
\end{equation}
otherwise we have
\begin{multline}
\label{eq:detD2}
\det\left(D_1(-e;N,l,e)\right)=(-1)^{e\cdot N}2^{\frac{2+e+e^2-2N}{2}}\po{e}{e}
	\prod_{j=1}^{N-e}(2j-1)!\times\\
	\prod_{j=1}^{e}(k-1)!^2(N+k-2e-1)!\po{\tfrac{k+e-1}{2}}{k-1}\po{k-N}{N-e}.
\end{multline}

In order to establish this we observe that $D_1(-e;N,l,e)$ has a block form
which is sketched in Figure~\ref{fig:matrix}. The figure has to be read
according to the following convention: If a block is bounded from
above by a horizontal line numbered $r_1$ at the left margin, is
bounded from below by a horizontal line numbered $r_2$ at the right
margin, is bounded from the left by a vertical line numbered $c_1$ 
at the bottom margin, and is bounded from the right by a vertical
line numbered $c_2$ at the top margin,
then the block consists of the entries from rows $i=r_1,\dots,r_2$ and
columns $j=c_1,\dots,c_2$.
This block form is easily established by routine verification
directly from the definitions \eqref{eq:def-D2-1} and \eqref{eq:def-D2-2}.

\begin{figure}
	\begin{picture}(150,190)(-150,-20)
		\put(0,0){\line(0,1){150}}
		\put(0,150){\line(1,0){150}}
		\put(150,150){\line(0,-1){150}}
		\put(150,0){\line(-1,0){150}}
		\put(50,0){\line(0,1){100}}
		\put(50,145){\line(0,1){5}}
		\put(100,0){\line(0,1){150}}
		\put(0,50){\line(1,0){100}}
		\put(145,50){\line(1,0){5}}
		\put(0,100){\line(1,0){150}}
		\put(0,100){\line(1,-1){100}}
		\put(122,45){{\Large 0}}
		\put(22,20){{\Large 0}}
		\put(30,80){{\Large 0}}
		\put(65,15){{\Large 0}}
		\put(71,71){{\Huge $*$}}
		\put(46,121){{\Huge $*$}}
		\put(120,120){{\Large $M$}}
		\put(13,63){{\Large $Q_2$}}
		\put(79,31){{\Large $Q_1$}}
		\put(2,-7){\tiny $1$}
		\put(52,-7){\tiny $N+1-2e$}
		\put(102,-7){\tiny $N+1-e$}
		\put(28,155){\tiny $N-2e$}
		\put(83,155){\tiny $N-e$}
		\put(145,155){\tiny $N$}
		\put(-5,145){\tiny $1$}
		\put(-19,95){\tiny $e+1$}
		\put(-36,45){\tiny $N+1-e$}
		\put(155,103){\tiny $e$}
		\put(155,53){\tiny $N-e$}
		\put(155,3){\tiny $N$}
	\end{picture}
\caption{}
\label{fig:matrix}
\end{figure}

Hence the determinant $D_1(m;N,l)$ factors as follows,
\begin{equation}
\label{eq:D2-factors}
\det\left(D_1(m;N,l)\right)=(-1)^{e(N-e)}\det(Q_2)\det(Q_1)\det(M).
\end{equation}

Since $Q_1$ and $Q_2$ are upper and lower triangular matrices,
respectively, 
it is easy to evaluate $\det(Q_1)$ and $\det(Q_2)$. We simply have to
multiply all the entries on the main diagonal. We obtain
\begin{equation}
\label{eq:evalQ2}
\det(Q_1)
=	\begin{cases}
		(-1)^{\frac{e(e-1)}{2}}\prod_{k=0}^{e-1}k!(2e-2k-1)!(N+k-2e)!
			&\text{if }l<N+1-e,\\
		0 &\text{if }l\geq N+1-e.
	\end{cases}
\end{equation}
Note that \eqref{eq:detD2zero} follows immediately from \eqref{eq:D2-factors}  
and \eqref{eq:evalQ2}.

By assumption we have $l\geq\frac{N+1}{2}>e$, by assertion \eqref{eq:detD2zero}
(which is already proved) we may assume $l\leq N-e$, so we are sure to
encounter row $l$ in $Q_2$. Multiplying the entries on the main diagonal gives
\begin{equation}
\label{eq:evalQ1}
\det(Q_2)=\left(\frac{1}{2}\right)^{N-2e+1}\prod_{j=1}^{N-2e}(j-1)!\po{N-2e-j+1}{N-j+1}.
\end{equation}
(Note that for $e=0$, $Q_2$ is the whole, unmodified matrix $D(0;N,l)$.)

For the evaluation of $\det(M)$ we employ Krattenthaler's Lemma (see
\cite[Lemma~2.2]{KratAM}),
\begin{multline}
\label{lem:Kratt}
\det_{1\le i,j\le n}\Big((X_i+A_n)\cdots(X_i+A_{j+1})
(X_i+B_j)\cdots (X_i+B_2)\Big)\\
\hskip2cm =\prod _{1\le i<j\le n} ^{}(X_i-X_j)\prod _{2\le i\le j\le n}
^{}(B_i-A_j),
\end{multline}
where $X_1,\dots,X_n$, $A_2,\dots,A_n$, and $B_2,
\dots B_n$ are arbitrary indeterminates.

The $(i,j)$-entry of matrix $M$ is the  
$(i,N-e+j)$-entry of matrix $D(-e;N,l)$. By assumption we have $l\geq\frac{N+1}{2}>e$, so we do not encounter row $l$
in $M$. Hence, we must evaluate
\begin{equation}
\det_{1\leq i,j\leq e}\left(
\tfrac{(-j-e+1)}{2}\po{N+i-j+1}{N-e+j-1}\po{2N-e-2i+j+2}{e-j}
\right).
\end{equation}
By taking out $\po{i-N}{N-e}$ from the $i$-th row and
$(-2)^{e-j}\frac{1-e-j}{2}$ from the $j$-th column, we get
\begin{multline}
\label{eq:M-Krattform}
\det_{1\leq i,j\leq e}\Biggl(
\left(i+(-1-N+\tfrac{1}{2})\right)\cdots\left(i+(-N-1+\tfrac{e-j}{2})\right)\times\\
\left(i+(-N-j+1)\right)\cdots\left(i+(-N-1)\right)
\Biggr),
\end{multline}
which is precisely of the form required for \eqref{lem:Kratt} (set $n=e$,
$X_i=i$, $A_j=-1-N+\frac{e-j+1}{2}$ and $B_j=-j-N+1$). So we obtain
\begin{equation}
\prod_{j=2}^{e}(j-1)!\po{\tfrac{e-j+1}{2}}{j-1}
\end{equation}
for the $(e\times e)$-determinant in \eqref{eq:M-Krattform}.  
Multiplying this with the factors pulled out previously, we have
\begin{equation}
\label{eq:evalM}
\det(M)=
(-2)^{\frac{e(e-1)}{2}-e}\po{e}{e}\prod_{i=1}^{e}
	(i-1)!\po{\tfrac{e-i+1}{2}}{i-1}\po{i-N}{N-e}.
\end{equation}

By inserting \eqref{eq:evalQ2}, \eqref{eq:evalQ1} and \eqref{eq:evalM} into
\eqref{eq:D2-factors} and simplifying the expression, we obtain
\eqref{eq:detD2}. Thanks to \eqref{eq:polydet}, this establishes
\eqref{eq:polyeval}, and thus completes Step~4.

\medskip
This finishes the proof of Lemma~\ref{lem:complexpart}.
\hfill\qed

\section{Open problems and conjectures}
In this paper, we computed the number
of rhombus tilings of a hexagon with sides $N,M,N,N,M,N$, which
contain a fixed rhombus on the symmetry axis. There are two questions
which suggest themselves:

\def\Item "#1"{\par\noindent\hangindent2\parindent%
  \hangafter1\setbox0\hbox{\rm#1\enspace}\ifdim\wd0>2\parindent%
  \box0\else\hbox to 2\parindent{\rm#1\hfil}\fi\ignorespaces}

\medskip
{\leftskip.5cm 

\Item "Question 1:" What is the number of rhombus tilings of a
hexagon with {\em arbitrary} side lengths $a,b,c,a,b,c$ which contain
an {\em arbitrary} fixed rhombus?

\Item "Question 2:" What is the number of rhombus tilings of a
hexagon with side lengths $a,b,c,a,b,c$ which contain
{\em several\/} fixed rhombi?

}
\medskip

Whereas it is too much to expect ``nice" answers to these questions in
general, there is indeed hope for further nice results in special
cases. We should mention that, on an asymptotic level, both questions
are settled by the Cohn, Larsen and Propp result 
\cite[Theorem~1]{CoLPAA} that was mentioned in the Introduction.

Regarding Question~1: Clearly, the approach that we used in this
paper is not good enough for any generalizations in these directions.
For, the use of the Matchings Factorization Theorem at the very
beginning requires a reflective symmetry of the region that we are
considering. So, under this approach we need to have a hexagon with
sides $N,M,N,N,M,N$ with the fixed rhombus on the symmetry axis which 
cuts through the sides of length $M$. (The reader should observe that
the Matchings Factorization Theorem cannot be used with respect to
the ``other" symmetry axis, i.e., the symmetry axis which runs
{\em parallel\/} to the sides of length $M$. For, this symmetry axis does not
have the required property that it ``separates" the dual graph, cf\@.
the statement of \cite[Theorem~1.2]{CiucAB}.)

However, it seems that the approach taken by Helfgott and Gessel
\cite{HeGeAA} does allow to obtain further results in this direction.
(Recall that they also obtained Theorems~\ref{thm:MEven} and
\ref{thm:MOdd} for the cases where the fixed rhombus is the central
rhombus.)
In fact, a further result, using the approach by Helfgott and Gessel,
has already been obtained by the authors and will be
subject of a forthcoming article \cite{FuKrAD}. There, we compute the
number of rhombus tilings of a hexagon with sides $N,M,N,N,M,N$, the
parameters $N$
and $M$ being of the {\em same} parity,
which contain a rhombus which is ``next to the center of the region".
(Note, that, since $N$ and $M$ have the same parity, there is no
central rhombus.)

\smallskip
Regarding Question~2, there are strong indications that, at least,
there are ``nice" results analogous to those of
Theorems~\ref{thm:MEven} and \ref{thm:MOdd} if we consider rhombus
tilings of a hexagon with sides $N,M,N,N,M,N$ which contain a given
{\em set\/} of rhombi on the symmetry axis which 
cuts through the sides of length $M$. 

Suppose that we fix rhombi
$l_1,l_2,\dots,l_r$ on the symmetry axis, and let
$L=\{l_1,l_2,\dots,l_r\}$. If we want to know the number of all
rhombus tilings which contain these rhombi, then we can use the same
approach as we used for proving Theorems~\ref{thm:MEven} and
\ref{thm:MOdd}. 

That is (see Section~2), we apply first the Matchings
Factorization Theorem. It implies that 
{\em the number
of rhombus tilings of a hexagon with sides $N,2m,N,N,2m,N$ which
contain the rhombi from the set $L$ equals
\begin{equation}\label{eq:MatchFactL}
2^{N-r}R(S(N-1,m))\tilde R(C(N,m,L)),
\end{equation}}\noindent
where $S(N-1,m)$ is the same ``reduced" ``upper region" as in Section~2, and
where $C(N,m,L)$ denotes the resulting
``lower half", i.e., a region similar to $C(N,m,l)$ (which appeared in
Figure~\ref{fig:2}), however, where along the ``upper border" of the
region the rhombi from the set $L$ are removed.
As in Section~2, the symbol $\tilde R(C(N,m,L))$
denotes the weighted count of rhombus tilings of $C(N,m,L)$ in which
each of the top-most (horizontal) rhombi counts with weight $1/2$.
A similar result (generalizing \eqref{eq:CiucuMOdd}) holds for the
case of a hexagon with sides $N,2m-1,N,N,2m-1,N$.

Clearly, since $R(S(N-1,m))$ is already known (see \eqref{eq:SNm1} and
\eqref{eq:SNm2}), ``all" we need is the weighted count $\tilde
R(C(N,m,L))$. Again, the tiling problem can be translated into
nonintersecting lattice paths, and from the lattice path
interpretation we obtain a determinant for $\tilde
R(C(N,m,L))$. To be precise, we have
\begin{equation}
\label{eq:CNmL1}
\tilde R(C(N,m,L))	=
	\det_{1\leq i,j\leq N}\left(
		\begin{cases}
\frac{(N+m-i)!}{(m+i-j)!\,(N+j-2i+1)!}(m+\frac{N-j+1}{2})&\text{ if
}i\notin L\\
			\frac{(N+m-i)!}{(m+i-j)!\,(N+j-2i)!}&\text{ if
}i\in L
		\end{cases}
	\right).
\end{equation}
The reader should compare this determinant with the one in
\eqref{eq:CNml1}. Now we do not have just {\em one} ``exceptional"
row, now we have $r=\vert L\vert $ ``exceptional" rows. In fact, this
determinant can be regarded as the ``mixture" of two matrices,
$$A=\left(\frac{(N+m-i)!}{(m+i-j)!\,(N+j-2i+1)!}(m+\tfrac{N-j+1}{2})
\right)
_{1\le i,j\le N}$$
and 
$$B=\left(\frac{(N+m-i)!}{(m+i-j)!\,(N+j-2i)!}\right)
_{1\le i,j\le N}.$$

Now, to perform the task of evaluating the determinant, we would try
to follow the proof of Lemma~\ref{lem:complexpart} in Section~3.
Indeed, Auxiliary Facts~I and II, and Steps~1--3 carry over, when
suitably modified. In particular, the analogue of \eqref{eq:polydef}
reads: 

{\em The determinant in \eqref{eq:CNmL1} equals
\begin{equation} \label{eq:polydefL}
\prod _{i=1}^{\floor{N/2}}
		\po{m+i}{N-2i+1}
\prod _{i=r} ^{\floor{N/2}}\po{m+i+\tfrac{1}{2}}{N-2i}\cdot
P(m;N,L),
\end{equation}
where $P(m;N,L)$ is a polynomial in $m$ of degree at most $r(N-r)$.}

The reader should compare with \eqref{eq:polydef} and observe the
differences: First, the second product in \eqref{eq:polydefL} is
``smaller" than the corresponding product in \eqref{eq:polydef}.
Second, and unfortunately, the degree of $P(m;N,L)$ is in general
considerably larger than the degree of $P(m;N,l)$. 

It is here where the problems start. Now it comes to carry over
Step~4. At present, we do not know how to accomplish that. In order
to determine $P(m;N,L)$ we would need $r(N-r)+1$ evaluations of
$P(m;N,L)$. What we are able to do is to follow Step~4 in Section~3
and determine the value of $P(m;N,L)$ at $m=0,-1,\dots,-N$.
Unfortunately, for $r>1$, this is not good enough.

Still, we do believe that a reasonable formula for $P(m;N,L)$ can be
found. In fact, it appears that, quite often,
$P(m;N,L)$ does indeed factor further. 
At this point, we want to direct the reader's attention
to the fact that this is also the case for $P(m;N,l)$, the polynomial
that was defined in \eqref{eq:polydef} and which equals the second line
of \eqref{eq:main}. For, from the expression given by the second line of 
\eqref{eq:main}, it is immediately obvious that $(m+l)_{N-2l+1}$
divides $P(m;N,l)$. In particular, if $l=1$ then $P(m;N,l)$ factors
completely into linear factors. 

A similar phenomenon seems to hold for $P(m;N,L)$ in general. Namely,
for $L=\{1,2,\dots,r\}$ it appears that $P(m;N,L)$ factors completely
into linear factors (see Conjecture~\ref{conj1} below). 
And, whenever we move a rhombus by 1 to the
right then the degree of $P(m;N,L)$ in $m$ increases by 2 (at least for
large $N$; compare Conjectures~\ref{conj1}--\ref{conj3} below). 

We have worked out a few conjectures corresponding to ``small" $L$,
meaning that the numbers in $L$ are small, which means that the fixed
rhombi are far left in the hexagon. 

\begin{conj}
\label{conj1}
Let $m$ be a nonnegative integer and $N$ and $r$ be positive integers
with $N\ge r$. The number
of rhombus tilings of a hexagon with sides $N,2m,N,N,2m,N$, which
contain rhombi $1,2,\dots,r$ on the symmetry axis 
which cuts through the sides of length $2m$, equals
\begin{multline}
2^{\frac {1} {2}(r-1)(r-2N)}\frac{\binom{m+N-1}{m}^2}{\binom{2m+2N-1}{2m}}
\prod _{i=N-r} ^{N-2}\frac {1} {i!}
\prod _{i=1} ^{r-1}\frac {(2i)!!\,(2N-2i-1)!!\,
(m+i+1)_{N-2i-1}} {(2i-1)!!\,(m+i+1/2)_{N-2i}}
\\\times
\prod_{i=1}^N\prod_{j=1}^N\prod_{k=1}^{2m}\frac{i+j+k-1}{i+j+k-2}.
\end{multline}
The number
of rhombus tilings of a hexagon with sides $N+1,2m-1,N+1,N+1,2m-1,N+1$, which
contain rhombi $1,2,\dots,r$ on the symmetry axis 
which cuts through the sides of length $2m-1$, equals
\begin{multline}
2^{\frac {1} {2}(r-1)(r-2N)}\frac{\binom{m+N-1}{m}^2}{\binom{2m+2N-1}{2m}}
\prod _{i=N-r} ^{N-2}\frac {1} {i!}
\prod _{i=1} ^{r-1}\frac {(2i)!!\,(2N-2i-1)!!\,
(m+i+1)_{N-2i-1}} {(2i-1)!!\,(m+i+1/2)_{N-2i}}
\\\times
\prod_{i=1}^{N+1}\prod_{j=1}^{N+1}\prod_{k=1}^{2m-1}\frac{i+j+k-1}{i+j+k-2},
\end{multline}
\end{conj}

\begin{conj}
\label{conj2}
Let $m$ be a nonnegative integer and $N$ and $r$ be positive integers
with $N\ge r+1$. The number
of rhombus tilings of a hexagon with sides $N,2m,N,N,2m,N$, which
contain rhombi $1,2,\dots,r-1,r+1$ on the symmetry axis 
which cuts through the sides of length $2m$, equals
\begin{multline}
2^{\frac {1} {2}(r-1)(r-2N)}\frac {3r(N-r)}
{(2r-1)(2N-2r+1)(m+r)(m+N-r)}\\
\times
\left(m^2+Nm+\tfrac {(2r-1)(2N-2r+1)} {3}\right)
\frac{\binom{m+N-1}{m}^2}{\binom{2m+2N-1}{2m}}
\prod _{i=N-r} ^{N-2}\frac {1} {i!}\\
\times
\prod _{i=1} ^{r-1}\frac {(2i)!!\,(2N-2i-1)!!\,
(m+i+1)_{N-2i-1}} {(2i-1)!!\,(m+i+1/2)_{N-2i}}
\prod_{i=1}^N\prod_{j=1}^N\prod_{k=1}^{2m}\frac{i+j+k-1}{i+j+k-2}.
\end{multline}
The number
of rhombus tilings of a hexagon with sides $N+1,2m-1,N+1,N+1,2m-1,N+1$, which
contain rhombi $1,2,\dots,r-1,r+1$ on the symmetry axis 
which cuts through the sides of length $2m-1$, equals
\begin{multline}
2^{\frac {1} {2}(r-1)(r-2N)}\frac {3r(N-r)}
{(2r-1)(2N-2r+1)(m+r)(m+N-r)}\\
\times
\left(m^2+Nm+\tfrac {(2r-1)(2N-2r+1)} {3}\right)
\frac{\binom{m+N-1}{m}^2}{\binom{2m+2N-1}{2m}}
\prod _{i=N-r} ^{N-2}\frac {1} {i!}\\
\times
\prod _{i=1} ^{r-1}\frac {(2i)!!\,(2N-2i-1)!!\,
(m+i+1)_{N-2i-1}} {(2i-1)!!\,(m+i+1/2)_{N-2i}}
\prod_{i=1}^{N+1}\prod_{j=1}^{N+1}\prod_{k=1}^{2m-1}\frac{i+j+k-1}{i+j+k-2},
\end{multline}
\end{conj}

\begin{conj}
\label{conj3}
Let $m$ be a nonnegative integer and $N$ and $r$ be positive integers
with $N\ge r+2$. The number
of rhombus tilings of a hexagon with sides $N,2m,N,N,2m,N$, which
contain rhombi $1,2,\dots,r-1,r+2$ on the symmetry axis 
which cuts through the sides of length $2m$, equals
\begin{multline}
2^{\frac {1} {2}(r-1)(r-2N)}\frac {45} {64}
\frac {(r)_2\,(N-r-1)_2} {(r-1/2)_2\,(N-r-1/2)_2\,
(m+r)_2\,(m+N-r-1)_2}\\
\left(m^4+2Nm^3+\left(N^2+\tfrac{(20r+1)N}9 -\tfrac {20r^2+2r+5}
{9}\right)m^2 \right.\kern5cm\\
\kern1cm\left.+\tfrac {\left((20r+1)N-20r^2-2r-5\right)} {9}Nm
+\frac {4} {45}(2r-1)(2r+1)(2N-2r-1)(2N-2r+1)\right)\\
\times
\frac{\binom{m+N-1}{m}^2}{\binom{2m+2N-1}{2m}}
\prod _{i=N-r} ^{N-2}\frac {1} {i!}
\prod _{i=1} ^{r-1}\frac {(2i)!!\,(2N-2i-1)!!\,
(m+i+1)_{N-2i-1}} {(2i-1)!!\,(m+i+1/2)_{N-2i}}
\\\times
\prod_{i=1}^N\prod_{j=1}^N\prod_{k=1}^{2m}\frac{i+j+k-1}{i+j+k-2}.
\end{multline}
The number
of rhombus tilings of a hexagon with sides $N+1,2m-1,N+1,N+1,2m-1,N+1$, which
contain rhombi $1,2,\dots,r-1,r+2$ on the symmetry axis 
which cuts through the sides of length $2m-1$, equals
\begin{multline}
2^{\frac {1} {2}(r-1)(r-2N)}\frac {45} {64}
\frac {(r)_2\,(N-r-1)_2} {(r-1/2)_2\,(N-r-1/2)_2\,
(m+r)_2\,(m+N-r-1)_2}\\
\left(m^4+2Nm^3+\left(N^2+\tfrac{(20r+1)N}9 -\tfrac {20r^2+2r+5}
{9}\right)m^2 \right.\kern5cm\\
\kern1cm\left.+\tfrac {\left((20r+1)N-20r^2-2r-5\right)} {9}Nm
+\frac {4} {45}(2r-1)(2r+1)(2N-2r-1)(2N-2r+1)\right)\\
\times
\frac{\binom{m+N-1}{m}^2}{\binom{2m+2N-1}{2m}}
\prod _{i=N-r} ^{N-2}\frac {1} {i!}
\prod _{i=1} ^{r-1}\frac {(2i)!!\,(2N-2i-1)!!\,
(m+i+1)_{N-2i-1}} {(2i-1)!!\,(m+i+1/2)_{N-2i}}
\\\times
\prod_{i=1}^{N+1}\prod_{j=1}^{N+1}\prod_{k=1}^{2m-1}\frac{i+j+k-1}{i+j+k-2},
\end{multline}
\end{conj}

Without difficulty we could move on and work out
further conjectures. Already these three conjectures do contain so
many similarities, so that the ``next" formula, the formula for fixing
rhombi $1,2,\dots,r-1,r+4$, can, almost, be guessed right away.
Still, these guesses do not seem to help for the
determinant evaluations that would be needed to prove these
enumerations. So we content ourselves with these three conjectures,
and hope that they give enough evidence that in this area there are a
lot of further beautiful results waiting to be unearthed, and
proved, of course.

\ifx\undefined\bysame
\newcommand{\bysame}{\leavevmode\hbox to3em{\hrulefill}\,}
\fi


\begin{thebibliography}{10}

\bibitem{CiucAB}
M.~Ciucu, {\em Enumeration of perfect matchings in graphs with reflective  
symmetry\/}, J. Combin\@. Theory Ser.~A {\bf 77} (1997), 67-97.

\bibitem{CiucKratAB}
M.~Ciucu and C.~Krattenthaler, {\em The number of centered lozenge tilings of a
symmetric hexagon\/}, preprint.

\bibitem{CoLPAA} 
H.    Cohn, M. Larsen and J. Propp, 
{\em The shape of a typical boxed plane partition},
preprint.

\bibitem{DT}
G. David and C. Tomei, {\em The problem of the calissons}, 
Amer\@. Math\@. Monthly\@. {\bf 96} (1989), 429--431.

\bibitem{FuKrAD}
M. Fulmek and C. Krattenthaler, {\em The number of rhombus
tilings of a symmetric hexagon which contain a fixed rhombus on the
symmetry axis, II}, in preparation.

\bibitem{GeViAB}
I.~M.~Gessel and X.~Viennot, {\em Determinants, paths, and plane  
partitions\/}, preprint, 1989.

\bibitem{HeGeAA} 
H.    Helfgott and I. M. Gessel,
{\em Exact enumeration of certain tilings of diamonds and hexagons with
defects}, preprint.

\bibitem{KratAM}
C.~Krattenthaler, {\em Generating functions for plane partitions of a given  
shape\/},  Manuscripta Math.\ {\bf 69}, (1990), 173--202.

\bibitem{KratAK} 
C.    Krattenthaler, {\em A determinant evaluation and some enumeration 
results for plane partitions}, in: Number-Theoretic Analysis,
E.~Hlawka, R.~F.~Tichy, eds., Lect\@. Notes in Math. 1452, 
Sprin\-ger-Ver\-lag, Berlin, 1990.

\bibitem{KratBG}
C.~Krattenthaler,
{\em Some $q$-analogues of determinant identities which arose in plane  
partition enumeration\/}, S\'eminaire Lotharingien Combin.\ {\bf 36},
(1996), paper~B36e, 23~pp.

\bibitem{KratBH}
C.~Krattenthaler, {\em A new proof of the M--R--R conjecture --- including a  
generalization\/},
preprint.

\bibitem{KratBI}
C.~Krattenthaler, {\em An alternative evaluation of the Andrews--Burge  
determinant\/}, to appear in the ``Ro\-ta\-fest\-schrift".

\bibitem{KratBD}
C.~Krattenthaler, {\em Determinant identities and a generalization of the  
number of totally symmetric self-complementary plane partitions\/},
Electron.\ J.\ Combin.\ {\bf 4}(1) (1997), \#R27, 62~pp.

\bibitem{KrZeAA}
C.~Krattenthaler and D.~Zeilberger, {\em Proof of a
determinant evaluation conjectured by Bombieri, Hunt and van der Poorten\/},  
New York J. Math., {\bf 3} (1997), 54--102.

\bibitem{MacMahon}
P.A.~MacMahon, {\em Combinatory Analysis\/}, vol.~2, Cambridge University
Press, 1916; reprinted by Chelsea, New York, 1960.

\bibitem{PropAA}
J.~Propp, {\em Twenty open problems on enumeration of matchings\/},
manuscript, (1996).

\bibitem{Prudnikov}
A.P.~Prudnikov, Yu.A.~Brychkov, O.I.~Marichev, {\em Integrals and Series\/},
vol.~3: Mores Special Functions, Gordon and Breach, New York, London, 1989.

\bibitem{SlatAC}
L.~J.~Slater,
{\em Generalized hypergeometric functions\/},
Cambridge University Press, Cambridge, (1996).

\bibitem{StemAE}
J.~R.~Stembridge, {\em Nonintersecting paths, pfaffians and plane partitions\/}, 
Adv.\ in Math.\ {\bf 83} (1990)
96---131.

\end{thebibliography}
\end{document}